\documentclass[12pt,hidelinks]{amsart}
\usepackage{bbm}
\usepackage{amsmath}
\usepackage{amsthm}
\usepackage{stmaryrd}
\usepackage{amsfonts,mathrsfs,amssymb}
\usepackage{verbatim}
\usepackage{tikz-cd}
\usepackage{times}
\usepackage{fullpage}
\usepackage{enumerate}
\makeatletter
\def\@tocline#1#2#3#4#5#6#7{\relax
  \ifnum #1>\c@tocdepth 
  \else
    \par \addpenalty\@secpenalty\addvspace{#2}%
    \begingroup \hyphenpenalty\@M
    \@ifempty{#4}{%
      \@tempdima\csname r@tocindent\number#1\endcsname\relax
    }{%
      \@tempdima#4\relax
    }%
    \parindent\z@ \leftskip#3\relax \advance\leftskip\@tempdima\relax
    \rightskip\@pnumwidth plus4em \parfillskip-\@pnumwidth
    #5\leavevmode\hskip-\@tempdima
      \ifcase #1
       \or\or \hskip 1em \or \hskip 2em \else \hskip 3em \fi%
      #6\nobreak\relax
    \hfill\hbox to\@pnumwidth{\@tocpagenum{#7}}\par
    \nobreak
    \endgroup
  \fi}
\makeatother

\title{Human Variations on AI Mathematics}
\bigskip

\author[J. Cao]{Junyan Cao}
\address{Laboratoire de Mathématiques J.A. Dieudonné UMR 7351 CNRS, Université Côte d'Azur Parc Valrose 06108, Nice, France}
\email{junyan.cao@unice.fr}
\urladdr{https://sites.google.com/site/junyancao}

\author[M. P\u aun]{Mihai P\u aun}
\address{Universit\"at Bayreuth, Mathematisches Institut, Lehrstuhl Mathematik VIII, Universit\"atsstrasse 30,
	D-95447, Bayreuth, Germany}
	\email{mihai.paun@uni-bayreuth.de}
	\urladdr{https://www.mathe8.uni-bayreuth.de/de/team/prof-paun/index.php}

\newcommand{\cE}{{\mathcal E}}

\newcommand{\cO}{{\mathcal O}}

\newcommand{\db}{\bar\partial}

\newcommand{\wt}{\widetilde}

\newcommand{\bbC}{{\mathbb C}}

\newcommand{\bbN}{{\mathbb N}}

\newcommand{\bbQ}{{\mathbb Q}}
\newcommand{\bbR}{{\mathbb R}}

\newcommand{\supp}{\operatorname{Supp}}

\newcommand{\Div}{\operatorname{Div}}

\newcommand{\dbar}{\bar \partial}
\newcommand{\ddbar}{\sqrt{-1}\partial\bar \partial}
\newcommand{\ep}{\varepsilon}

\newcommand{\ord}{{\rm ord}}

\def\XXint#1#2#3{{\setbox0=\hbox{$#1{#2#3}{\int}$}
\vcenter{\hbox{$#2#3$}}\kern-.5\wd0}}

\usepackage{hyperref}

\newcommand{\step}[2]{\par\medskip\noindent\textbf{Step #1: #2.}\ }

\begin{document}

\theoremstyle{plain}
\newtheorem{m-thm}{\sc Theorem}
\newtheorem*{s-thm}{\sc Theorem}
\newtheorem*{s-cor}{\sc Corollary}
\newtheorem{lem}{\sc Lemma}[section]
\newtheorem{thm}[lem]{\sc Theorem}
\newtheorem{prop}[lem]{\sc Proposition}
\newtheorem{cor}[lem]{\sc Corollary}
\newtheorem{claim}[lem]{\sc Claim}

\theoremstyle{definition}
\newtheorem{conj}[lem]{\sc Conjecture}
\newtheorem*{sconj}{\sc Conjecture}
\newtheorem{prob}[lem]{\sc Open Problem}
\newtheorem{defn}[lem]{\sc Definition}
\newtheorem*{s-defn}{\sc Definition}
\newtheorem{qn}[lem]{\sc Question}
\newtheorem{ex}[lem]{\sc Example}
\newtheorem{rmk}[lem]{\sc Remark}
\newtheorem*{s-ex}{\sc Example}
\newtheorem*{s-qn}{\sc Question}
\newtheorem*{s-prop}{\sc Proposition}
\newtheorem*{s-rmk}{\sc Remark}
\newtheorem*{s-rmks}{\sc Remarks}
\newtheorem{rmks}[lem]{\sc Remarks}
\newtheorem*{ack}{\sc Acknowledgments}
\newtheorem*{notation}{\sc \underline{Notation}}
\newcommand{\RR}{\mathbb R}

\renewcommand\thesubsubsection{\thesubsection.\roman{subsubsection}}

\begin{abstract} This article is motivated by the OpenAI manuscripts \emph{Minimal
metrics and interior injectivity for nef adjoints} and \emph{Log
abundance in characteristic zero}. We establish a jet comparison for
metrics with minimal singularities and use it to prove a dlt-type
extension theorem. Combined with geometric arguments from the latter manuscript, this yields an alternative analytic proof of its
Theorem~3.1.
\end{abstract}

\maketitle

\tableofcontents

\section{Introduction}

The main result of this article is a divisorially log terminal (dlt)-type extension theorem
for pluricanonical sections. We show that under two positivity assumptions on the
log canonical divisor, the restriction map to the reduced boundary
on a log resolution is surjective. As an application, we give an
alternative proof of Theorem~3.1 of the OpenAI manuscript
\cite{OpenAI2}, \emph{Log abundance in characteristic zero}:
our extension result replaces the formal lifting and algebraization
steps in the original argument.
\smallskip

The starting point was a reading of the recently released OpenAI
manuscripts \cite{OpenAI1,OpenAI2}, prompted by the question of how
their arguments relate to analytic methods in complex geometry.
The proof developed here brings ideas from these manuscripts
together with the methods of Ohsawa--Takegoshi,
Demailly--Hacon--P\u{a}un, and Cao--Demailly--Matsumura.
It relies on two analytic results: a comparison theorem for singular
metrics, yielding the vanishing of the Lelong numbers of a metric
with minimal singularities, and an extension theorem from a reduced
dlt boundary.
\smallskip

\noindent
Our main extension theorem is the following.

\begin{thm}[Theorem~\ref{extensionthm}]\label{extensionthmintro}
Let $(Y,\Delta_Y)$ be a projective $\mathbb Q$-factorial dlt pair. Assume that $K_Y+\Delta_Y$ is nef and that
$K_Y+(1-\delta_0)\Delta_Y$ is pseudo-effective for some $\delta_0>0$.
Let $\pi:X\to Y$ be a log resolution, and write
\[
K_X+S+B=\pi^*(K_Y+\Delta_Y)+E,
\]
where $S$ is reduced, $B,E\geq0$ are rational divisors,
$\lfloor B\rfloor=0$, $E$ is $\pi$-exceptional, and no two of
$S,B,E$ share an irreducible component. Assume that their union has
simple normal crossings, and set
\[
F:=\pi^*(K_Y+\Delta_Y).
\]
Then, for every sufficiently divisible integer $m\geq2$, the
restriction map
\[
H^0(X,mF)\longrightarrow H^0(S,mF|_S)
\]
is surjective. Equivalently, every section of
$m(K_X+S+B)|_S$ in the image of multiplication by $s_{mE}|_S$
extends to $X$.
\end{thm}

Before discussing the proof, we turn to the application. Theorem~3.1 of \cite{OpenAI2} establishes abundance for nef log canonical divisors admitting a
\emph{signed representative} (\emph{machine dixit}) supported on
the boundary, subject to the additional hypotheses below.
It is a central step in the proof of the main theorem of that
manuscript. Here is the statement.

\begin{thm}[\cite{OpenAI2}, Theorem~3.1]\label{extenopenai}
Let $(Y,\Delta_Y)$ be a projective $\mathbb Q$-factorial dlt pair
with reduced boundary
\[
\Delta_Y=\sum_j D_j,
\]
and set $L:=K_Y+\Delta_Y$. Assume that $L$ is nef, that
$K_Y+(1-\delta_0)\Delta_Y$ is pseudo-effective for some
$\delta_0>0$, and that
\[
L\sim_{\mathbb Q}\sum_j a_jD_j=:G,
\qquad a_j\in\mathbb Q.
\]
If $L|_{\Delta_Y}$ is semiample, then
\[
\kappa(Y,L)=\nu(Y,L).
\]
\end{thm}

The deduction from our extension theorem can be outlined as follows;
a detailed proof is given in Corollary~\ref{coropenAI}. We may
assume that
\[
0<\nu(Y,L)<\dim Y.
\]
By \cite[Lemma~3.2]{OpenAI2}, the numerical dimension drops by
at least one upon restriction to any boundary component, and by
exactly one on at least one component occurring in $G$ with positive
coefficient. Moreover, the morphism defined by a globally generated
multiple of $L|_{\Delta_Y}$ maps the intersection of the positive
and nonpositive parts of the boundary to a set of dimension strictly
less than $\nu(Y,L)-1$.

This allows us to construct sections on the boundary that vanish
on its nonpositive part and define an image of dimension
$\nu(Y,L)-1$ on a positive component. Our extension theorem lifts
these sections to $Y$. The signed representative, after cancellation
of its poles, provides a further global section vanishing on that
component. Together, these sections raise the image dimension by
one, giving
\[
\kappa(Y,L)\geq\nu(Y,L).
\]
\medskip

\noindent
We now describe the analytic ingredients in the proof of the
extension theorem, Theorem~\ref{extensionthmintro}.

The first is the vanishing of all Lelong numbers of a positively
curved metric $\varphi_{F,\mathrm{min}}$ with minimal singularities
on $F$; see Theorem~\ref{nu0}. Its proof uses a key idea taken from
\cite{OpenAI1}, \emph{Minimal metrics and interior injectivity
for nef adjoints}. For a klt pair $(X,\Delta)$ with
$K_X+\Delta$ nef, that manuscript proves that a positively curved
metric with minimal singularities on $K_X+\Delta$ has zero Lelong
numbers everywhere; see Definition~\ref{defnminsing} for the
terminology. Pursuing the method further, we obtain a jet version
that dispenses with nefness; this is developed in
Section~\ref{jetversion}.

The central estimate is the following. Let $(X,B)$ be a klt pair
with $K_X+B$ pseudo-effective, fix a sufficiently ample line bundle
$A$, and let $\varphi_{\min}$ be a positively curved metric with
minimal singularities on $K_X+B$. We prove that
\begin{equation}\label{jetest}
\nu(\varphi_{\min},x_0)
\leq
\liminf_m
\frac1m
\min_{0\neq s\in V_m}
\operatorname{ord}_{x_0}(s)
\end{equation}
at every point $x_0\in X$, where
\[
V_m:=H^0(X,m(K_X+B)+A).
\]

Passing to limits of normalized section weights does not suffice
to establish \eqref{jetest}. Although such a limit defines a
semipositively curved metric on $K_X+B$, its Lelong numbers may
increase: upper semicontinuity yields an inequality opposite to
the one required. The adjoint structure of $K_X+B$ is what makes
the sharper control possible.

The proof combines an extremal-section argument with an $L^2$
estimate, in the spirit of the Ohsawa--Takegoshi extension theorem
\cite{OT}. The cutoff argument in \cite{OpenAI1}, however, has a
distinctive feature. In the Ohsawa--Takegoshi proof, the right-hand
side of the estimate converges to the $L^2$-norm on the hypersurface;
here it converges to zero. A nontrivial section is nevertheless
recovered by contradiction.

The argument is also closely related to the proof of the
$L^{2/m}$ extension theorem. A section supplies the singular weight
for an $L^2$ estimate, and H\"older's inequality turns the resulting
bound into control of the pluricanonical integral. This interplay
between sections and singular metrics has long been central to
extension arguments; see, for example, \cite{Siu02,Pau07,DHP}.

We apply \eqref{jetest} to $K_Y+(1-\varepsilon)\Delta_Y$ and
let $\varepsilon\to0$. Lemma~\ref{small} supplies the required
bound on normalized vanishing orders, uniformly in $\varepsilon$.
We conclude that
\begin{equation}\label{minsing}
\nu(\varphi_{F,\mathrm{min}},x_0)=0
\qquad\text{for every }x_0\in X,
\end{equation}
which proves Theorem~\ref{nu0}.
\medskip

The second ingredient is an extension theorem
combining the methods of \cite{DHP} and \cite{CDM}. For a
sufficiently divisible integer $m\geq2$, set
\[
M:=(m-1)(K_X+S+B)+B,
\]
and equip this line bundle with the metric $h_M$ induced by
$h_{F,\mathrm{min}}$ and the canonical singular metrics of $B$
and $E$. Then
\[
m(K_X+S+B)=K_X+S+M.
\]
We seek to extend a section
\[
s\in H^0(S,(K_X+S+M)|_S)
\]
divisible by $s_E^{\otimes m}|_S$.

Following a key idea of \cite{DHP}, we combine $h_{F,\mathrm{min}}$
with the pseudo-effectivity of $K_Y+(1-\delta_0)\Delta_Y$ to
construct a metric $h_S$ on $\mathcal O_X(S)$ satisfying
\[
\sqrt{-1}\Theta_{h_M}(M)\geq0,
\qquad
\sqrt{-1}\Theta_{h_M}(M)\geq c_0\,\sqrt{-1}\Theta_{h_S}(S)
\]
for some $c_0>0$. Although the weight of $h_S$ need not be
quasi-psh, the minimality of the singularities of
$\varphi_{F,\mathrm{min}}$ ensures that $|s_S|_{h_S}$ is uniformly
bounded above. We thus obtain the required curvature inequalities
and norm bound.

To treat the snc divisor $S$ and the singular metrics, we adapt the
argument by approximate solutions in \cite{CDM}.
Simultaneous regularization preserves the curvature inequalities,
while the vanishing of the Lelong numbers supplies the local $L^2$
control. These estimates yield the desired extension to $X$.

\smallskip

\noindent \textbf{Organization of the article.}
In Section~\ref{secpre}, we recall the notion of minimal singularities and an $L^2$ estimate used throughout the paper.
In Section~\ref{secminI}, we prove a local generation criterion for the vanishing of Lelong numbers of metrics with minimal singularities, following the method of \cite{OpenAI1}.
In Section~\ref{jetversion}, we generalize this result to the pseudo-effective setting.
Finally, in Section~\ref{mainextension}, we prove our main extension theorem, Theorem~\ref{extensionthmintro}, and give a complete proof of its application to Theorem~3.1 of \cite{OpenAI2}.
\smallskip

\noindent \textbf{AI disclosure.} The analytic content of this article was developed and written by us. ChatGPT assisted with language and exposition, with explanations of the proof of Theorem 3.1 in \cite{OpenAI2}, including its connection with our main result on the extension of pluricanonical forms. 

\section{Preliminaries}\label{secpre}

Let $X$ be a compact K\"ahler manifold, and let $\alpha$ be a smooth real $d$-closed $(1,1)$-form on $X$. We assume that the cohomology class  $\{\alpha\}\in H^{1,1}(X, \bbR)$ is pseudo-effective, meaning that it contains a closed positive current. The following notion plays an important role in this paper.
\begin{defn}\label{defnminsing}
The psh upper envelope $V_\alpha$ of the class $\{\alpha\}$ is defined as
\[V_\alpha:=\left(\sup\{u\in\operatorname{PSH}(X,\alpha):u\leq0\}\right)^\star.\]
Here the star denotes upper-semicontinuous regularization. Then
$\alpha+\ddbar V_\alpha$ is a positive current with minimal
singularities.
\end{defn}

\noindent If $\alpha=\sqrt{-1}\Theta(L,h_0)$ for a smooth metric $h_0$,
the metric $h_{\min}=h_0e^{-V_\alpha}$ has minimal singularities
among the positively curved metrics on $L$. Its local weights will be denoted by $\varphi_{\rm min}$.
\smallskip

\noindent We recall a few properties and applications of this notion.
\begin{enumerate}[\rm (i)]
\smallskip

\item It plays a key role in establishing extension theorems for pluricanonical forms, starting with the work of Siu.
\smallskip

\item If the class $\{\alpha\}$ is big, then the Laplacian of
$V_\alpha$ is bounded by a function depending on the potential of a K\"ahler current with analytic singularities $T\in \{\alpha\}$, by a result of \cite{BD}.
\smallskip

\item If the class $\{\alpha\}$ is nef and big, $V_\alpha$ can be obtained via a family of Monge--Amp\`ere equations (and this plays an important role in the recent work of Guedj, Guenancia, Guo, Phong, Song, and Sturm, among others).
\end{enumerate}
\smallskip

We first introduce the cutoff functions used in the following sections, then recall a version of H\"ormander's $L^2$ estimates, due to Demailly, suited to our applications.

\subsection{Cutoff functions}\label{cutoff} Fix a constant $C_\infty>0$ and a smooth cutoff function $\chi:\bbR\to[0,1]$ satisfying
\begin{equation}\label{eq20}
\chi|_{]-\infty, -1-2C_\infty]}= 1,\qquad  \chi|_{[-2C_\infty, \infty[}= 0.
\end{equation}
We also introduce an auxiliary convex function.

For $c>0$, choose a smooth, convex, increasing function
$\rho_c:\bbR\to\bbR$ with the following properties:
\begin{enumerate}[\rm (i)]

\item \emph{It is affine with slope $c$ on the interval} $]-\infty, -2-2C_\infty]$.
\smallskip

\item \emph{It is strictly convex on the interval} $]-1-2C_\infty, -2C_\infty]$.
\smallskip

\item \emph{It is affine on the interval} $]-2C_\infty+1,\infty[$.
\end{enumerate}
\smallskip

Such a function $\rho_c$ exists and satisfies
\begin{equation}\label{eq21}
\rho_c(t)\leq C_1+ M_ct_+, \qquad t_+:= \max (t, 0)
\end{equation}
for some constant $C_1>0$, where
\[M_c:=\sup\rho_c'<\infty.\]
\smallskip

\noindent For a prescribed $c>0$, the function $\rho_c$ has
strictly positive second derivative on a
neighborhood of $\supp\chi'$.
Thus,
\[\sup_{\supp(\chi')}\frac{|\chi'|^2}{\rho_c''}<\infty.\]
Throughout the applications below, $c$ and $\rho_c$ are fixed first,
and then the positive integer $m$ will be chosen so that $m-1>M_c$.
\medskip

\subsection{\texorpdfstring{$L^2$}{L2} estimates}\label{L2} We recall the following standard result from \cite{Dem12}.

\begin{thm}\label{hjp} Let $(\Omega,\omega)$ be an $n$-dimensional Kähler manifold that also admits a complete Kähler metric. Let $(L,h)$ be a holomorphic line bundle with a smooth Hermitian metric satisfying $\sqrt{-1}\Theta(L,h)\geq0$, and denote by $A_1$ the operator $[\sqrt{-1}\Theta_h(L),\Lambda_\omega]$.
Let $v$ be an $L$-valued $(n,1)$-form satisfying
$$v\in L^2, \qquad \langle A^{-1}_1v, v\rangle \in L^1, \qquad \dbar v=0.$$
Here the $L^p$ spaces are defined using the Hermitian structures
$(\Omega, \omega)$ and $(L, h)$.

\noindent Then there exists an $L$-valued $(n,0)$-form $u$ such that
\[\dbar u= v, \qquad \int_\Omega |u|^2e^{-\varphi_L}\leq
\int_\Omega \langle A^{-1}_1v, v\rangle_\omega e^{-\varphi_L} dV_\omega.\]
\end{thm}
The inverse quadratic form above is understood to be infinite on the kernel of $A_1$; thus its integrability includes the corresponding range condition. The singular metric case will be treated by approximation in Section~\ref{regsubsection}.
\smallskip

\noindent We will apply this statement in the following context. Let $X$ be a projective manifold.  Let $B$ be a holomorphic $\bbQ$-line bundle endowed with a singular metric $h_B =e^{-\varphi_B}$ with positive curvature current. Let $m\geq2$ be an integer such that $mB$ is a line bundle,
and consider a holomorphic section
\begin{equation}\label{eq70}
s\in H^0(X, m(K_X+ B)+ A),\qquad \int_X|s|^{\frac{2}{m}}e^{-\varphi_B- \frac{\varphi_A}{m}}< \infty
\end{equation}
with finite $L^{2/m}$ integral, where $A$ is a line bundle equipped with a smooth, semipositively curved metric $e^{-\varphi_A}$. In particular, this includes an ample bundle with a positive metric and the pullback of such a bundle under a modification. The argument below uses only semipositivity of this metric.

Assume that $K_X+B$ is pseudo-effective, and fix a positively curved singular metric $h=e^{-\phi}$ on it. We regard $s$ as an $(n,0)$-form with values in the line bundle
\[L_m:= B+ \frac{m-1}{m}\bigl(m(K_X+B)+A\bigr)+ \frac{1}{m}A\]
which we equip with the metric of local weight
\begin{equation}\label{eq71}
\phi_m:= \varphi_B+  \frac{m-1}{m}\log|f|^2 + \frac{1}{m}\varphi_A .
\end{equation}
Here $f$ denotes the local coefficient of $s$.
\smallskip

Consider the $L^1$ function $\displaystyle \Psi:= \phi+\frac{1}{m}\varphi_A- \frac{1}{m}\log |f|^2$, as well as the modified weight
\begin{equation}\label{eq72}
\phi_m+ \rho_c\circ \Psi
\end{equation}
where $\rho_c$ is the function introduced in Section~\ref{cutoff}. Denote the resulting metric on $L_m$ by $h_m$.

\noindent We claim that
\begin{equation}\label{eq73}
\sqrt{-1}\Theta(L_m, h_m)\geq \sqrt{-1}\rho_c''(\Psi)\partial\Psi\wedge \dbar\Psi
\end{equation}
holds weakly on $X$. First, \eqref{eq73} holds away from the zero divisor of $s$, since
$$\ddbar \log |f|^2= 0 \qquad\text{on }X\setminus \rm Div (s) .$$
On the other hand,
\begin{equation}\label{eq74}
\rho_c\circ \Psi\leq C+ M_c\max(0, \phi+\frac{1}{m}\varphi_A- \frac{1}{m}\log |f|^2)\leq C+M_c\Big(-\frac{1}{m}\log|f|^2\Big)
\end{equation}
after shrinking the chart so that $|f|\leq1$ and adjusting $C$.

On a relatively compact trivializing chart, the upper bound
\[
\phi_m+\rho_c(\Psi)\leq C+\frac{m-1-M_c}{m}\log|f|^2
\]
follows from \eqref{eq74} and the local upper bound for the fixed psh weight. This proves local upper boundedness and psh extension across $f=0$. The singular right-hand side of \eqref{eq73} can be understood by the regularization argument below; no undefined product of distributions is needed.

\smallskip

The right-hand side of our $\dbar$ equation is the $L_m$-valued $(n,1)$-form
\[v:=\dbar(\chi\circ\Psi)\wedge s.\]
Using the notation of Theorem~\ref{hjp}, a direct calculation gives
\begin{equation}\label{eq75}
\int_X\langle A_1^{-1}v, v\rangle_\omega e^{-\phi_m- \rho_c\circ \Psi}dV_\omega\leq C\int_{S_\Psi}|s|^{\frac{2}{m}}
e^{-\varphi_B- \frac{\varphi_A}{m}}<\infty
\end{equation}
since the integral on the left is supported on $\supp v$, which is contained in
\begin{equation}\label{eq84}
S_\Psi:= \big(-1- 2C_\infty\leq \Psi\leq - 2C_\infty\big).
\end{equation}
\smallskip

Although $\dbar v=0$, the singularities of $h_m$ prevent a direct application of Theorem~\ref{hjp}. The standard regularization procedure applies in both the projective and Kähler settings. We briefly recall it in the projective case.
\subsubsection{Regularization procedure}\label{regsubsection} Consider
an effective ample divisor $H$ of sufficiently high degree containing the divisor $(s=0)$ as well as the divisors needed to trivialize the rational bundles. Choose
\begin{equation}\label{eq76}
X\setminus H= \bigcup_k \Omega_k, \qquad \overline{\Omega}_k\subset \Omega_{k+1}
\end{equation}
as an exhaustion of $X\setminus H$ by relatively compact, smoothly bounded, strictly pseudoconvex domains.
On these domains, we can write our weights
\begin{equation}\label{eq77}
\phi|_{\Omega_k}= \lim \phi_\ep, \qquad \varphi_B= \lim \varphi_{B, \ep}
\end{equation}
as decreasing limits of smooth, psh functions by smooth psh approximation on the chosen Stein domains. By the properties of $H$, all the bundles involved, namely $K_X$, $B$, and $A$, are trivial when restricted to $X\setminus H$. On $\Omega_k$ we then introduce
\begin{equation}\label{eq78}
\Psi_\ep:= \phi_\ep+\frac{1}{m}\varphi_A- \frac{1}{m}\log |f|^2, \qquad v_\ep:= \dbar (\chi\circ\Psi_\ep)\wedge s
\end{equation}
as well as the metric $h_{m, \ep}$ given by the weight
\begin{equation}\label{eq79}
\phi_{m, \ep}+ \rho_c\circ \Psi_\ep, \qquad \phi_{m, \ep}:= \varphi_{B, \ep}+ \frac{1}{m}\varphi_A+ \frac{m-1}{m}\log|f|^2.
\end{equation}
Now all the weights and forms in \eqref{eq78} and \eqref{eq79} are smooth,
and the fundamental estimate
\begin{equation}\label{eq80}
\int_{\Omega_k}\langle A_1^{-1}v_\ep, v_\ep\rangle_\omega e^{-\phi_{m, \ep}- \rho_c\circ \Psi_\ep}dV_\omega\leq C\int_{S_{\Psi_\ep}\cap \Omega_k}|s|^{\frac{2}{m}}
e^{-\varphi_B- \frac{\varphi_A}{m}}
\end{equation}
still holds, with right-hand side bounded uniformly in $k$ and $\varepsilon$.

Theorem~\ref{hjp} gives a solution of the $\dbar$ equation satisfying the uniform estimates
\begin{equation}\label{eq81}
\dbar u_{\ep,k}= v_\ep, \qquad \int_{\Omega_k}|u_{\ep,k}|^2e^{-\phi_{m, \ep}- \rho_c\circ \Psi_\ep}\leq C\int_{S_{\Psi_\ep}\cap \Omega_k}|s|^{\frac{2}{m}}
e^{-\varphi_B- \frac{\varphi_A}{m}}.
\end{equation}
We first let $\ep\to0$. We claim that the limit satisfies
\begin{equation}\label{eq82}
\dbar u_k= v, \qquad \int_{\Omega_k}|u_k|^2e^{-\phi_{m}- \rho_c\circ \Psi}\leq C\int_{S_{\Psi}\cap \Omega_k}|s|^{\frac{2}{m}}
e^{-\varphi_B- \frac{\varphi_A}{m}}
\end{equation}
for every $k$. Here both the $\dbar$ equation and the integration domain must be considered in passing to the limit.

Indeed, the regularized weights are
\[
\varphi_{B,\varepsilon}+\frac{\varphi_A}{m}
+\frac{m-1}{m}\log|f|^2+\rho_c(\Psi_\varepsilon),
\]
and decrease as $\ep\to0$. The left-hand side of the inequality in \eqref{eq82} therefore passes to the limit by weak lower semicontinuity.

Next, the functions $\chi(\Psi_\varepsilon)$ are bounded and converge almost everywhere to $\chi(\Psi)$. Hence
\[
\chi(\Psi_\varepsilon)s\longrightarrow\chi(\Psi)s
\quad\hbox{locally in }L^1.
\]
Their $\bar\partial$-derivatives in the sense of distributions converge. For a locally weakly convergent subsequence of $(u_{\ep,k})_\ep$, the limit therefore satisfies
\[
\bar\partial u_k=\bar\partial\bigl(\chi(\Psi)s\bigr).
\]
This formulation \emph{defines} the right-hand side without requiring an ordinary derivative of the singular weight.

Finally, we choose $\chi'$ supported in a compact interval contained in $(-1-2C_\infty,-2C_\infty)$. Almost-everywhere convergence and domination by the fixed integrable density bound the limiting right-hand side by its integral over a closed annulus contained in $\{\Psi<-2C_\infty\}$.
\smallskip

\noindent Letting $k\to\infty$, we obtain a solution $u$ of the original equation satisfying
\begin{equation}\label{eq83}
\dbar u= v, \qquad \int_X|u|^2e^{-\phi_{m}- \rho_c\circ \Psi}\leq C\int_{S_{\Psi}}|s|^{\frac{2}{m}}
e^{-\varphi_B- \frac{\varphi_A}{m}}
\end{equation}
Notice that the modified section $\chi(\Psi)s- u$
is locally $L^2$ with respect to a smooth
metric and extends holomorphically across the
divisor we removed. This also proves the global distributional equation.

\section{Metrics with minimal singularities, I}\label{secminI}

Let $X$ be a projective manifold, and let $\Delta$ be a klt $\mathbb Q$-divisor on $X$.
Let $\varphi_\Delta$ be the weight of the canonical singular metric on $\Delta$, namely $\ddbar \varphi_\Delta =[\Delta]$.  Assume also that $A\gg 0$ is a sufficiently positive ample line bundle on $X$. We show that the following local version of \cite[Thm 1.1]{OpenAI1} holds.

\begin{thm}\label{numin}
Let $x_0 \in X$ be a point such that, for every sufficiently large and divisible integer $m$ clearing the denominators of $\Delta$, there exists a section $s_m \in H^0 (X, m(K_X+\Delta)+A)$ that does not vanish at $x_0$.
Then the metric with minimal singularities $h_{\rm min} = e^{-\varphi_{\rm min}}$
on $K_X+\Delta$ has zero Lelong number at $x_0$.
\end{thm}
\noindent Under the hypothesis of this theorem, the adjoint divisor $K_X+\Delta$ is automatically pseudo-effective.

\begin{proof}
We follow \cite{OpenAI1}, clarifying a few points and references along the way. We fix a smooth metric of weight $\varphi_A$ on $A$ such that $\ddbar \varphi_A >0$ on $X$.
On the finite-dimensional vector space
\[V_m:= H^0\big(X, m(K_X+ \Delta)+ A\big)\]
we define the length
\begin{equation}\label{eq11}
\Vert \cdot\Vert_m: V_m\to \bbR, \qquad \Vert s\Vert_m:= \left(\int_X |s|^{2/m}e^{-\varphi_\Delta- \frac{\varphi_A}{m}}\right)^{m/2},
\end{equation}
where $\varphi_\Delta$ is the weight of the canonical singular metric on $\Delta$.
For $m\geq2$, $\Vert\cdot\Vert_m$ need not be a norm, but it is the natural length function associated with the $L^{2/m}$ Bergman kernel; see \cite{BP08}. A standard compactness argument shows that, for each $m$, there exists a section $s_m\in V_m$ such that
\begin{equation}\label{eq12}
 |s_m(x_0)|= \max_{\Vert s\Vert_m= 1} |s(x_0)|,
\end{equation}
so that $s_m$ computes the value of the Bergman kernel mentioned above.
\medskip

\noindent Suppose, for a contradiction, that
\begin{equation}\label{eq13}
\nu(\varphi_{\rm min}, x_0)>0.
\end{equation}
We will construct, for some $m\gg0$, a section $\sigma_m\in V_m$ satisfying
\begin{equation}\label{eq14}
\sigma_m(x_0)= s_m(x_0), \qquad \Vert \sigma_m\Vert_m< 1
\end{equation}
contradicting the extremality of $s_m$.
\medskip

\noindent First, after passing to a subsequence, we obtain a limit
\begin{equation}\label{eq15}
\varphi_\infty:= \lim_m\frac{1}{m}\log |f_m|^2
\end{equation}
which defines the local weight of a metric on $K_X+\Delta$. Here $f_m$ is the local coefficient of $s_m$.
Indeed, the normalization and the mean-value inequality give local uniform upper bounds for the section weights. We may therefore extract a subsequence converging in $L^1$ and almost everywhere. The same upper bounds give an integrable majorant, so dominated convergence yields
\begin{equation}\label{eq16}
\int_X e^{\varphi_\infty -\varphi_\Delta}= 1,
\end{equation}
which excludes the alternative in which the weights tend uniformly to $-\infty$ on compact subsets.
By the minimality of the singularities of $\varphi_{\rm min}$, there exists a constant $C_\infty>0$ such that
\begin{equation}\label{eq17}
\varphi_{\rm min} + C_{\infty}\geq \varphi_\infty
\end{equation}
holds on $X$.
\medskip

We now prepare to apply the Ohsawa--Takegoshi-type argument of \cite[Lemma~2.3]{OpenAI1}. We introduce the global function
\begin{equation}\label{eq18}
\Psi_m:= -\frac{1}{m}\log\bigl(|f_m|^2e^{-m\varphi_{\rm min}-\varphi_A}\bigr)
\end{equation}
on $X$. By \eqref{eq17}, the quantity
\begin{equation}\label{eq19}
\ep_m:= \int_{\Psi_m< -2C_\infty}|s_m|^{2/m}e^{-\varphi_\Delta- \frac{\varphi_A}{m}}
\end{equation}
tends to zero as $m\to\infty$. This follows from dominated convergence, since \eqref{eq15} and \eqref{eq17} hold almost everywhere.

We now construct $\sigma_m$ satisfying \eqref{eq14}. Write
\[m(K_X+ \Delta)+ A=: K_X+ L_m\]
and endow $L_m$ with the metric whose local weights are
\begin{equation}\label{eq22}
\phi_m:= \varphi_\Delta+ \frac{m-1}{m}\log |f_m|^2+ \frac{1}{m}\varphi_A+ \rho_c\circ \Psi_m
\end{equation}
where $c> 0$ is a positive real number such that
\begin{equation}\label{eq86}
e^{-c\varphi_{\rm min}}\not\in L^1(X, x_0).
\end{equation}
This is where we use the fact that the Lelong number of $\varphi_{\rm min}$
at $x_0$ is strictly positive.
As we have already seen in \eqref{eq74}, the weights $\phi_m$ are locally bounded from above
and moreover
\begin{equation}\label{eq23}
\ddbar \phi_m\geq \sqrt{-1}\rho_c''(\Psi_m)\partial \Psi_m\wedge \dbar\Psi_m.
\end{equation}

\medskip
As in the proof of Ohsawa--Takegoshi \cite{OT}, we consider the following $L_m$-valued $(n,1)$-form
\begin{equation}\label{eq24}
v_m:= \dbar(\chi_m s_m)\qquad \chi_m:= \chi\circ \Psi_m.
\end{equation}
Since $s_m$ is holomorphic, this can be rewritten as
\begin{equation}\label{eq25}
v_m= \chi'(\Psi_m)\dbar\Psi_m\wedge s_m.
\end{equation}
The results of Section~\ref{L2} give an $L_m$-valued form $u_m$ satisfying
\begin{equation}\label{eq26}
\dbar u_m= v_m, \qquad \int_X|u_m|^2e^{-\phi_m}\leq \int_X|v_m|^2_{\ddbar \phi_m}e^{-\phi_m}
\end{equation}
where the subscript means that we use the form $\ddbar \phi_m$ to measure the antiholomorphic part of $v_m$.
\smallskip

Now we estimate $\int_X|u_m|^2e^{-\phi_m}$.
Since the support of $v_m$ is contained in the set
\begin{equation}
S_{\Psi_m}:= -1- 2C_\infty\leq {\Psi_m}\leq -2C_\infty
\end{equation}
by the choice of $\chi$, equations \eqref{eq83} and \eqref{eq19} give
\begin{equation}\label{eq28}
\int_X|u_m|^2e^{-\phi_m}\leq C\ep_m
\end{equation}
where $C$ is a positive constant, independent of $m$.

Since $s_m(x_0)\ne0$ and $\nu(\varphi_{\min},x_0)>0$, we have $\chi_m=1$ on a neighborhood of $x_0$. In particular, $u_m$ is holomorphic there. The choice of $c$ in \eqref{eq86}, together with its weighted integrability, forces $u_m(x_0)=0$. Consequently the holomorphic section
\begin{equation}\label{eq29}
\sigma_m:= \chi_m s_m- u_m
\end{equation}
has the \emph{same} value as the extremal section $s_m$ at $x_0$.  Applying H\"older's inequality twice to \eqref{eq28} gives
$\Vert\sigma_m\Vert_m< 1$. We refer to Step 5 of the proof of Theorem~\ref{generamin} for more details on this point. Thus $\sigma_m$ satisfies \eqref{eq14}, contradicting extremality and proving the theorem.
\end{proof}

\section{Metrics with minimal singularities, II}\label{jetversion}

\noindent In this section, we establish a jet version of Theorem~\ref{numin}. The setting is as follows.
\begin{enumerate}[(1)]
    \item Let $Y$ be a normal projective variety and let $\pi: X\to Y$ be a desingularization. Let $B\geq0$ be an effective $\mathbb Q$-divisor on $X$ such that $(X,B)$ is klt and
\begin{equation*}
K_X +B \equiv_{\mathbb R} F + E,
\end{equation*}
where $E\geq0$ is rational, $\pi$-exceptional, and $F$ is the pullback of some $\mathbb Q$-Cartier divisor on $Y$. Let $A_Y$ be a sufficiently ample line bundle on $Y$.
\smallskip

\item Let $x_0\in X$ be an arbitrary point. For each sufficiently divisible positive integer $m$ such that $mF$ is a line bundle, set
\[V_m:= H^0\big(X, m F + \pi^* A_Y\big).\]
Define
\[\delta_m:= \frac{1}{m}\min_{s\neq 0}\{\ord_{x_0}(s) : s\in V_m\}\]
the minimum normalized vanishing order at $x_0$ among nonzero sections in $V_m$. If $V_m=0$, we set $\delta_m=\infty$.
\end{enumerate}

\medskip

\noindent The main result of this section is the following.
\begin{thm}\label{generamin}
In addition to (1) and (2) above, assume that $F$ is pseudo-effective
and let $\varphi_{F,\rm min}$ be the weight of a positively curved metric with minimal singularities on $F$. Then
$$\nu(\varphi_{F,\rm min},x_0)\leq\liminf_m\delta_m.$$
\end{thm}
\medskip

\noindent Before proving the theorem, we formulate a possible Kähler analogue.

\begin{qn}
Let $X$ be a compact Kähler manifold and let $B$ be a klt divisor on $X$. We suppose that $K_X +B$ is pseudo-effective.
Consider a smooth representative $\alpha\in c_1(K_X+B)$ and write
$T_{\min}=\alpha+\ddbar V_\alpha$.
Similarly, take a current with minimal singularities $T_{\rm min, \ep}\in c_1 (K_X +B) +\ep \omega_X$.
Is it true that the equality \[
\nu(T_{\min},x)=\lim_{\varepsilon\downarrow0}
\nu(T_{\min,\varepsilon},x)
\]
holds for every $x\in X$?
\end{qn}

\begin{proof}[Proof of Theorem~\ref{generamin}]
We adopt the following convention: we interpret every numerical bundle identification using a fixed
unitary flat metric on the numerically trivial rational bundle
$K_X+B-F-E$. This does not affect curvature. All multiples $m$ below are assumed to clear
the denominators of $B,E,F$ and of this flat twist. The notation
$A=\pi^*A_Y$ will be used.
\smallskip

\noindent We define the length function
	\begin{equation}\label{eq30}
		\Vert \cdot\Vert_m:  V_m\to \bbR, \qquad \Vert s\Vert_m:= \left(\int_X |s|^{2/m}|s_E|^2e^{-\varphi_{B}- \frac{\varphi_A}{m}}\right)^{m/2}
	\end{equation}
Fix local coordinates $z$ centered at $x_0$ and a local frame of the line bundle. Consider the jet evaluation function
	$$\cE(s)= \max_{|\alpha|= m\delta_m} \Big|\frac{\partial^\alpha (s)}{\partial z^\alpha}(x_0)\Big|$$
	defined on $V_m$. Whenever $V_m\neq0$, there exists a section $s_m\in V_m$ with $\Vert s_m\Vert_m= 1$ such that
	\begin{equation}\label{eq31}
		\cE(s_m)= \max_{\Vert s\Vert_m= 1} \cE(s).
	\end{equation}
	In other words, $s_m$ computes the value of the jet version of the $L^{2/m}$ Bergman kernel.
	This remains true despite the presence of $\displaystyle s_{E}$ in \eqref{eq30}.
	\medskip

Suppose, for a contradiction, that
	\begin{equation}\label{eq32}
		\nu(\varphi_{F,\rm min}, x_0)> \lim\inf_m \delta_m .
	\end{equation}
We will construct, for some $m\gg0$, a section $\sigma_m\in V_m$ such that
	\begin{equation}\label{eq33}
		\cE(\sigma_m)= \cE(s_m), \qquad \Vert \sigma_m\Vert_m< 1
	\end{equation}
	contradicting the choice of $s_m$ and proving the theorem.
    The proof is divided into several steps.
We first pass to a subsequence on which $\delta_m\to d<\nu(\varphi_{F,\rm min}, x_0)$, and choose the constants
$M,\eta>0$ with
\[\delta_m\leq M, \qquad \nu(\varphi_{F,\rm min}, x_0)-\delta_m\geq\eta.\]
All the limits in what follows are taken along this subsequence.
We choose the slope $c$ of the convex increasing function $\rho_c$ so that
\[
c\eta>M+\ord_{x_0}(E)+n,
\]
see Step 4.
	\medskip

\step{1}{Extracting a limit} As in the previous section, we can extract a limit
	\begin{equation}\label{eq34}
		\varphi_\infty:= \lim_m\big(\log|s_{E}|^2+\frac{1}{m}\log |s_m|^2\big)
	\end{equation}
	representing the weights of a metric on $K_X+ B$. Indeed, the normalization $\Vert s_m\Vert_m=1$ and the mean-value inequality give local uniform upper bounds for the weights on the right-hand side of \eqref{eq34}.
	Moreover, by Lebesgue dominated convergence we infer that
	\begin{equation}\label{eq35}
		\int_X e^{\varphi_\infty -\varphi_{B}}= 1,
	\end{equation}
	so in particular $\varphi_\infty$ cannot be identically $-\infty$.
	\smallskip

	The weights $\varphi_\infty$ define a positively curved metric on $K_X+B$. Since $E$ is $\pi$-exceptional and $F$ is pulled back from $Y$, we can write
	\begin{equation}\label{eq36}
		\varphi_\infty= \varphi_{F, \infty}+ \varphi_{E}
	\end{equation}
	where $\varphi_{E}$ is the weight associated to the canonical singular metric on $E$, and $\varphi_{F, \infty}$ is the weight of a metric on $F$ such that $\ddbar \varphi_{F, \infty}\geq 0$.
	\smallskip

	By minimality of $\varphi_{F, {\rm min}}$, it follows that there exists a positive constant $C_\infty> 0$ such that
	\begin{equation}\label{eq37}
		\varphi_{F, {\rm min}} + C_{\infty}\geq \varphi_{F, \infty}
	\end{equation}
	holds on $X$. Now we introduce the global function
	\begin{equation}\label{eq38}
		\Psi_m:= - \frac{1}{m}\log (|s_m|^2e^{-m\varphi_{F, {\rm min}}- \varphi_A})
	\end{equation}
	on $X$ (the factor corresponding to $E$ cancels against the weight $\varphi_E$). The quantity
	\begin{equation}\label{eq39}
		\ep_m:= \int_{\Psi_m< -2C_\infty}|s_m|^{2/m}|s_E|^2e^{-\varphi_{B}- \frac{\varphi_A}{m}}
	\end{equation}
	tends to zero as $m\to\infty$. This follows from dominated convergence, since \eqref{eq34} holds almost everywhere.
	\medskip

	\step{2}{$L^2$ estimates} We now construct $\sigma_m$ satisfying \eqref{eq33}. Recall that
 \[T:= K_X+ B-F-E\]
 is a numerically trivial $\bbQ$-bundle, which we endow with a flat metric.

    Write
	\[m(F+ E)+ A=: K_X+ L_m\]
	and notice that $s_m s_E^m$ becomes an $L_m= \big((m-1)(K_X+B)+B+A- mT\big)$-valued $(n,0)$-form. With respect to unitary flat frames for $T$, we define the following metric on $L_m$
	\begin{equation}\label{eq42}
		\phi_m:= \varphi_{B}+ \frac{m-1}{m}\log |f_m f_E^m|^2+ \frac{1}{m}\varphi_A+ \rho_c\circ \Psi_m.
	\end{equation}
The parameter $c$ is chosen as specified above.
As explained in the previous section, we have
	\begin{equation}\label{eq43}
		\ddbar \phi_m\geq \sqrt{-1}\rho_c ''(\Psi_m)\partial \Psi_m\wedge \dbar\Psi_m \qquad\text{ on } X.
	\end{equation}

	As in the proof of Ohsawa--Takegoshi \cite{OT}, we consider the following $L_m$-valued $(n,1)$-form
	\begin{equation}\label{eq44}
		v_m:= \dbar(\chi_m \cdot s_m\otimes s_{E}^{\otimes m}),\qquad \text{where }\chi_m:= \chi\circ \Psi_m.
	\end{equation}
	Since $s_m$ is holomorphic, this can be rewritten as
	\begin{equation}\label{eq45}
		v_m= \chi'(\Psi_m)\dbar\Psi_m\wedge \big(s_m\otimes s_{E}^{\otimes m}\big)
	\end{equation}
We claim that $\displaystyle |v_m|_{\ddbar \phi_m}^2e^{-\phi_m}\in L^1(X)$.

Indeed, this expression defines a global $(n,n)$-form on $X$. Since the support of $v_m$ is contained in $\Psi_m< -2C_\infty$, we have
\[|v_m|_{\ddbar \phi_m}^2e^{-\phi_m}\leq C|s_m\otimes s_E^{\otimes m}|^2e^{-\varphi_{B}- \frac{m-1}{m}\log |f_m f_E^m|^2-\frac{1}{m}\varphi_A}\]
on the support of $v_m$, thanks to the curvature lower bound and the
fact that the antiholomorphic part of $v_m$ is proportional to $\dbar\Psi_m$. The right-hand side equals
\[C|s_m|^{\frac{2}{m}} |s_E|^2e^{-\varphi_{B}-\frac{1}{m}\varphi_A}\]
which is clearly in $L^1$.
    \smallskip

The $L^2$ estimates of Section~\ref{L2} give an $L_m$-valued form $u_m$ satisfying
	\begin{equation}\label{eq46}
		\dbar u_m= v_m, \qquad \int_X|u_m|^2e^{-\phi_m}\leq \int_X|v_m|^2_{\ddbar \phi_m}e^{-\phi_m}
	\end{equation}
	where the subscript means that we use the form $\ddbar \phi_m$ to measure the antiholomorphic part of $v_m$. Concerning the $L^2$-norm of $u_m$, notice that the support of $v_m$ is contained in the set
	\begin{equation}\label{eq47}
		-1- 2C_\infty\leq {\Psi_m}\leq -2C_\infty
	\end{equation}
	by the choice of the cutoff function $\chi$. The estimate in \eqref{eq46} becomes 
	\begin{equation}\label{eq48}
	\int_X|u_m|^2e^{-\phi_m}\leq C\ep_m
	\end{equation}
	where $C$ is a positive constant, independent of $m$ -- this is done as in the previous section.
	\medskip

\step{3}{Construction of the competing section $\sigma_m$}
	\noindent Consider next the holomorphic section
	\begin{equation}\label{eq49}
		\widetilde\sigma_m:= \chi_m s_m s_{E}^m- u_m.
	\end{equation}
	Although this section may a priori vanish identically, as a holomorphic section of
	\[m(F+ E)+ \pi^\star A_Y\]
	it is divisible by $s_{E}^{\otimes m}$. Thus we can write
	\[u_m=: w_m\otimes s_{E}^m\]
	for a global measurable section $w_m$ of $mF+\pi^\star A_Y$, satisfying
	\begin{equation}\label{eq50}
		\int_{(X, x_0)}|w_m|^2|s_{E}|^2e^{-\frac{m-1}{m}\log |f_m|^2-\rho_c\circ \Psi_m -\varphi_B- 1/m\varphi_A}<\infty,
	\end{equation}
	as a consequence of \eqref{eq46}.
	We consider
	\begin{equation}\label{eq51}
		\sigma_m:= \chi_m s_m- w_m \in H^0 (X, m F +\pi^* A_Y).
	\end{equation}
	In the next two steps, we show that $\sigma_m$ satisfies \eqref{eq33} for $m$ large enough.
\medskip

\step{4}{Jets of $\sigma_m$}
We prove that
$$\cE(\sigma_m)=\cE(s_m).$$
By the change-of-variables formula, for almost every line $\ell=(tv)_{t\in\bbC}$ through $x_0=0$,
	\begin{equation}\label{eq64}
		\int_{(\bbC, 0)}F_m(tv)|t|^{2n-2}d\lambda(t)<\infty,
	\end{equation}
	by \eqref{eq50}, where
	\begin{equation}\label{eq65}
		F_m:= |w_m|^2|s_E|^2e^{-\frac{m-1}{m}\log |f_m|^2-\rho_c\circ \Psi_m -\varphi_B- 1/m\varphi_A}.
	\end{equation}
Moreover, when restricted on a
full-measure set of lines, the section $s_m$ has order $m\delta_m$. Set
\[\nu:= \nu(\varphi_{F,\min}, x_0),\quad e: =\ord_{x_0}(E), \quad b:=\ord_{x_0}(B)\]
and note that
\[
\varphi_{F,\min}|_\ell=\nu\log|t|^2+\tau_\ell
\]
with $\tau_\ell$ subharmonic and bounded above, by Siu's slicing theorem. Moreover, the aforementioned 
set can be chosen simultaneously for the countable subsequence of $m$'s.
\smallskip

Since $\nu-\delta_m\geq\eta>0$, the function $\Psi_m|_\ell$
tends to $-\infty$. Thus $\chi_m=1$ near zero on the line, and
$w_m|_\ell=(s_m-\sigma_m)|_\ell$ is holomorphic. If its order
$\ell_m$ were at most $m\delta_m$, the weighted integral would dominate
a constant multiple of
\[
\int_{|t|<r}|t|^{2\ell_m+2e-2(m-1)\delta_m
-2c(\nu-\delta_m)-2b+2n-2}\,d\lambda(t).
\]
This integral diverges, because
\[
\ell_m+e-(m-1)\delta_m-c(\nu-\delta_m)-b+n
\leq M+e+n-c\eta<0
\]
by the choice of parameters at the beginning of the proof.
\smallskip

\noindent Consequently $\ell_m> m\delta_m$ on almost every such line.

Finally, consider the difference
\[g_m:=s_m- \sigma_m\]
and notice that it is holomorphic.
 The homogeneous
Taylor polynomials of $g_m$ through degree $m\delta_m$ vanish in almost every direction and hence vanish identically. Thus $s_m$ and $\sigma_m$ have the same jet through this order, and
$\cE(\sigma_m)=\cE(s_m)$.

\begin{rmk}
The introduction of $g_m$ is essential: we cannot argue directly with $w_m$, which is a priori only measurable. On a generic line, however, $g_m$ coincides with $w_m$ near the origin, allowing us to conclude.
\end{rmk}

	\medskip

\step{5}{Norm of $\sigma_m$}  Finally, we show that
$\Vert\sigma_m\Vert_m< 1$ for $m$ large enough.
We define
\[
I_m(s):=\int_X|s|^{2/m}|s_E|^2e^{-\varphi_B-\varphi_A/m},
\qquad \|s\|_m=I_m(s)^{m/2},
\]
so that $I_m(s_m)=1$. In compatible local frames, set
\[
d\mu_m:=|s_E|^2e^{-\varphi_B-\varphi_A/m},\qquad
\tau_m:=\max\left(\varphi_{F,\min}+\frac{\varphi_A}{m},
                      \frac1m\log|f_m|^2\right).
\]
The expressions below are interpreted invariantly, including the unitary flat metric when present. Recall that $M_c=\sup\rho_c'<\infty$, and assume $m>M_c+1$. Since
\[
\rho_c(\Psi_m)\leq C+M_c\left(\tau_m-\frac1m\log|f_m|^2\right),
\]
the $L^2$ estimate gives
\[
\int_X|w_m|^2|s_m|^{-2(m-1-M_c)/m}e^{-M_c\tau_m}\,d\mu_m
\leq C\varepsilon_m.
\]
On the support of $\chi_m$, the quantity
$\tau_m-m^{-1}\log|f_m|^2=\max(\Psi_m,0)$ equals zero. Thus the corresponding integral for $\chi_ms_m$ is bounded by $\varepsilon_m$. Since $\sigma_m=\chi_ms_m-w_m$, we obtain
\begin{equation}\label{eq:weightedchallenger}
\int_X|\sigma_m|^2|s_m|^{-2(m-1-M_c)/m}e^{-M_c\tau_m}\,d\mu_m
\leq C\varepsilon_m.
\end{equation}
Applying H\"older's inequality to \eqref{eq:weightedchallenger},
\[
I_m(\sigma_m)
\leq(C\varepsilon_m)^{1/m}J_m^{(m-1)/m},
\quad
J_m:=\int_X|s_m|^{2(m-1-M_c)/(m(m-1))}
                  e^{M_c\tau_m/(m-1)}\,d\mu_m.
\]
A second application, using $I_m(s_m)=1$, gives
\[
J_m\leq
\left(\int_X e^{\tau_m}\,d\mu_m\right)^{M_c/(m-1)}.
\]
The last integral is uniformly bounded. Indeed,
\[
e^{\tau_m}\leq e^{\varphi_{F,\min}+\varphi_A/m}+|s_m|^{2/m},
\]
and the first term is integrable against $d\mu_m$, uniformly with respect to $m$, because the fixed psh weight is locally bounded above, $E$ is effective, and $(X,B)$ is klt; the second has integral one. Consequently,
\[
I_m(\sigma_m)^m\leq C\varepsilon_m
\]
where $C$ is independent of $m$.
For large $m$, $I_m(\sigma_m)<1$, hence $\|\sigma_m\|_m<1$. This contradicts extremality and completes the proof.

\end{proof}

\section{Extension of pluricanonical forms}\label{mainextension}

\noindent We now apply the preceding results to prove an extension theorem for twisted pluricanonical forms. We work in the following setting.
\medskip

\noindent \textbf{Setting:} Let $Y$ be a normal projective variety, and consider the following data.

\begin{enumerate}[\rm (i)]

\item There is an effective $\mathbb Q$-divisor $\Delta_Y$ such that $(Y,\Delta_Y)$ is a $\bbQ$-factorial dlt pair.
\smallskip

\item The adjoint divisor $K_Y+\Delta_Y$ is nef, and there exists a constant $\delta \in ]0,1[$ such that
\[K_Y+ (1-\delta)\Delta_Y\]
is pseudo-effective.
\end{enumerate}
\medskip

\noindent Choose a log resolution $\pi: X\to Y$ of $(Y, \Delta_Y)$ and define
\begin{equation}\label{eq59}
K_X+ \Delta_X:= \pi^\star(K_Y+ \Delta_Y).
\end{equation}
Note that the support of $\Delta_X$ has simple normal crossings and its coefficients are all less than or equal to one. We decompose $\Delta_X$ as follows:
$$
\Delta_X := S +B -E,
$$
where $S$ consists of the components with coefficient one, $B$ is the klt part, and $-E$ is the negative part. Then we have
\begin{equation}\label{eq590}
K_X+ S +B = F +E ,
\end{equation}
where $F:= \pi^\star(K_Y+ \Delta_Y)$ is nef and $E$ is $\pi$-exceptional.

\medskip

\noindent We first prove the following result.

\begin{thm}\label{nu0} Assume that (i) and (ii) hold. Then the metric with minimal singularities $h_{\rm min}= e^{-\varphi_{\rm min}}$ on the $\bbQ$-bundle $K_X+ S+ B$ admits a decomposition
\[\varphi_{\rm min}= \varphi_F+ \varphi_E\]
where $\varphi_E$ is the weight of the canonical singular metric on $\cO(E)$. The metric $e^{-\varphi_{F}}$ on $F$ is semipositively curved and $\nu(\varphi_F, x_0)= 0$
for any $x_0\in X$.
\end{thm}

\begin{rmk}
Theorem~\ref{nu0} may be viewed as a log version of the main theorem of \cite{OpenAI1}. Since $E$ is $\pi$-exceptional, any positively curved metric $h$ on $K_X+S+B$ admits a decomposition of the form
\begin{equation}\label{canonicaldecom}
i\Theta_h(K_X+S+B)=T+[E],
\end{equation}
where $T\geq 0$ is a semipositive current in the class of $F$. We need to show that, if $h$ has minimal singularities, then the Lelong numbers of $T$ vanish everywhere.
\end{rmk}

\noindent Before proving Theorem~\ref{nu0}, we introduce some notation.
First, note that $Y$ is klt. Indeed, since $\Delta_Y$ is
an effective $\mathbb Q$-Cartier divisor, the pair $(Y,0)$ is still dlt.
Since a dlt pair with zero boundary is klt, it follows that $Y$ is klt;
see \cite[Proposition 2.41]{KM98}.
We define
$$
\Delta_X^{(0)}:=\pi^\star K_Y-K_X.
$$
Then every coefficient of $\Delta_X^{(0)}$ is strictly less than one.
Moreover, the difference
$$
\Delta_X-\Delta_X^{(0)}
=\pi^\star\Delta_Y\geq 0
$$
is effective.
For every $0<\varepsilon<\delta$, we set
\[
\Delta_{X,\varepsilon}
:=\Delta_X-\varepsilon\pi^\star\Delta_Y
=(1-\varepsilon)\Delta_X+\varepsilon\Delta_X^{(0)}.
\]
Then all coefficients of $\Delta_{X,\varepsilon}$ are strictly less
than one, so $(X,\Delta_{X,\varepsilon})$ is sub-klt.
Moreover, we write
\begin{equation}\label{eq60}
	F_\ep:=\pi^\star\bigl(K_Y+(1-\ep)\Delta_Y\bigr).
\end{equation}
Then $K_X+\Delta_{X,\ep}=F_\ep$ is pseudo-effective.

Next, we decompose $\Delta_{X,\varepsilon}$ into its positive and
negative parts:
\begin{equation}\label{eq62}
	B_\ep:=(\Delta_{X,\ep})_+,
	\qquad
	E_\ep:=(-\Delta_{X,\ep})_+.
\end{equation}
Thus,
\begin{equation}\label{eq58}
	K_X+B_\ep=F_\ep+E_\ep.
\end{equation}
Here, $E_\varepsilon$ is exceptional because the coefficients along
the strict transforms are $(1-\varepsilon)$ times the nonnegative
coefficients of $\Delta_Y$.
Notice that the family of divisors $(\Delta_{X,\ep})_{\ep>0}$
is monotone. Therefore, as $\ep$ decreases to zero, $(B_\ep)$
increases and $(E_\ep)$ decreases. Consequently, we have
$$
\lim_{\ep\to 0}B_\ep=S+B,
\qquad
\lim_{\ep\to 0}E_\ep=E,
\qquad
\lim_{\ep\to 0}F_\ep=F.
$$
\medskip

\noindent We first prove the following statement.
\begin{lem}\label{small}
For every fixed point $x_0\in X$, there exist an ample
divisor $A_Y$ on $Y$ and a constant $C> 0$ such that, for rational
$0<\ep\leq\delta$ and sufficiently divisible $m$, there
is a nonzero section
\[
s_{m,\ep}\in H^0(X,mF_\ep+\pi^*A_Y)
\]
such that
\[
\frac1m\ord_{x_0}(s_{m,\ep})
\leq C\left(\ep+\frac1m\right)
\]
holds. We stress that the constant $C$ is independent of both $\ep$ and $m$.
\end{lem}

\begin{proof} The lemma can be proved either algebraically or by $L^2$ methods. We give the algebraic argument here; the analytic argument is presented in the appendix.

Fix $x_0\in X$, and assume, after decreasing $\delta$ if necessary, that $\delta$ is rational. Set $F_Y:=K_Y+\Delta_Y$. By hypothesis, $F_Y$ is nef, and the preceding discussion shows that $Y$ is klt.
We can use Kawamata--Viehweg vanishing and Castelnuovo--Mumford regularity for $F_Y$ on $Y$. For instance, choose a very ample $H$ and an ample $A_Y$ so that $A_Y-K_Y-iH$ is ample for $1\leq i\leq n$. Since $p F_Y$ is nef, the required higher cohomology groups vanish for every admissible $p$.
It follows that
$pF_Y+A_Y$ is globally generated for every sufficiently divisible $p\geq0$.
We may therefore choose a section
\[\sigma_p\in H^0(X,pF+\pi^*A_Y), \qquad \sigma_p(x_0)\neq 0\]
nonvanishing at $x_0$.

The pseudo-effective divisor $F_\delta$ admits a positively curved singular
metric. Then there exists an ample bundle $A_X$ on $X$ such that $H^0(X,qF_\delta+A_X)\neq0$ for every sufficiently
divisible $q>0$, see, for example, the argument in the appendix. Moreover, by increasing $A_{Y}$ if necessary, we can assume that $\pi^*A_Y-A_X$ has a nonzero section. It follows that
\[
0\neq t_q\in H^0(X,qF_\delta+\pi^*A_Y).
\]
A degree bound for the corresponding effective zero divisor gives
\[\ord_{x_0}(t_q)\leq C_1q+C_2\]
with constants independent of $q$.

For each fixed rational $0<\varepsilon\leq\delta$, we take $m$ divisible enough that
\[
p=m(1-\varepsilon/\delta),\qquad q=m\varepsilon/\delta
\]
satisfy the required divisibilities. Now the section $\sigma_p$ is chosen nonvanishing at $x_0$. As explained above, we have a nonzero section
$t_q$ of $qF_\delta+\pi^*A_Y$. Then their product satisfies
\[
\sigma_pt_q\in H^0(X,mF_\varepsilon+2\pi^*A_Y),\qquad
m^{-1}\operatorname{ord}_{x_0}(\sigma_pt_q)
\leq (C_1/\delta)\varepsilon+C_2/m.
\]
Finally, replace $A_Y$ by $2A_Y$. The key point is that the constants are independent of $\varepsilon, m$, while divisibility may depend on the fixed rational $\varepsilon$. This proves the lemma.\end{proof}

\noindent We next turn to the proof of Theorem~\ref{nu0}.
\begin{proof}[Proof of Theorem~\ref{nu0}]
We use the notation introduced before Lemma~\ref{small}.
For each rational number $0<\varepsilon\leq\delta$, we have
\[
K_X+B_\varepsilon=F_\varepsilon+E_\varepsilon ,
\]
 where $E_\varepsilon$ is $\pi$-exceptional and $B_\varepsilon$ is klt.
By Theorem~\ref{generamin} and Lemma~\ref{small}, it follows that
\[
\nu(\varphi_{F_\varepsilon,\min},x_0)\leq C\varepsilon,
\]
where $C$ is independent of $\varepsilon$.

Let $\varphi_{F,\min}$ be the weight of a positively curved metric with minimal singularities on $F$. The positively curved metric with weight
$\varphi_{F_\varepsilon,\min}+\varepsilon\varphi_{\pi^\star\Delta_Y}$
on $F=F_\varepsilon+\varepsilon \pi^\star\Delta_Y$ and the minimality of $\varphi_{F,\min}$ imply that
\[
\varphi_{F,\min}\geq
\varphi_{F_\varepsilon,\min}+\varepsilon\varphi_{\pi^\star\Delta_Y}-C_\varepsilon
\]
for some constant $C_\ep$ which may depend on $\varepsilon$; however, this has no effect on Lelong numbers. Thus
\[
0\leq\nu(\varphi_{F,\min},x_0)
\leq\varepsilon\bigl(C+\ord_{x_0}\pi^\star\Delta_Y\bigr).
\]
Letting $\varepsilon\downarrow0$ proves that $\displaystyle \nu(\varphi_{F, {\rm min}}, x_0)=0$. Finally, note that the decomposition \eqref{canonicaldecom} identifies the minimal metric on $F+E$ with the sum of a minimal metric on $F$ and the canonical singular metric on $E$. The theorem is thus proved.
\end{proof}
\medskip

\noindent We now use Theorem~\ref{nu0} to prove the main extension theorem.
We recall the setting for convenience.

\begin{thm}\label{extensionthm}
Let $\pi\colon X\to Y$ be a log resolution of a projective
$\mathbb Q$-factorial dlt pair $(Y,\Delta_Y)$. Write
\[
 F_Y=K_Y+\Delta_Y,\qquad F=\pi^*F_Y,\qquad
 K_X+S+B=F+E,
\]
where $S$ is reduced, the coefficients of $B$ belong to $[0,1)$, and $E\geq0$ is $\pi$-exceptional. Assume that no two of $S,B,E$ share an irreducible component and that their union has simple normal crossings.
Assume also that:
\begin{enumerate}[\rm (1)]
\item $F_Y$ is nef;
\smallskip

\item $F_Y-\delta\Delta_Y$ is pseudo-effective
for some $\delta\in ]0,1[$.
\end{enumerate}
Then, for every sufficiently divisible integer $m\geq2$, the map
\[
 H^0(X,\cO_X(mF))\longrightarrow H^0(S,\cO_S(mF))
\]
is surjective. Equivalently, every section of
$m(K_X+S+B)|_S$ which is divisible by $s_E^m|_S$ extends to $X$.\end{thm}

\begin{proof} We divide the proof into six steps. We may assume that $S\neq0$.

\step{1}{Construction of metrics}
Let $\varphi_{F,\rm min}$ be a positively curved weight with minimal
singularities on $F$. Thanks to Theorem~\ref{nu0}, we have
\begin{equation}\label{nulelong}
\nu(\varphi_{F,\rm min},x)=0 \qquad\text{for every } x\in X.
\end{equation}

We claim that every component of $S$ is contained in the support of $\pi^*\Delta_Y$.
Indeed, write the pullbacks for the two pairs as
\[
K_X+S+B-E=\pi^\star(K_Y+\Delta_Y),\qquad
K_X+\Delta_X^{(0)}=\pi^\star K_Y.
\]
Since $(Y,0)$ is klt, the coefficients of $\Delta_X^{(0)}$ are strictly less than one. Each component $S_i$ of $S$ has coefficient one in $\Delta_X$, so
\[
\operatorname{coeff}_{S_i} \pi^*\Delta_Y
=1-\operatorname{coeff}_{S_i}\Delta_X^{(0)}>0.
\]
Consequently, there exists $\ep_0>0$ such that
\begin{equation}\label{eq90}
\pi^*\Delta_Y=\ep_0 S+R,\qquad R\geq0.
\end{equation}
Here $R$ is an effective $\mathbb R$-divisor.

Set
$$F_\delta:=\pi^* (F_Y-\delta \Delta_Y).
$$
As we assumed that $F_Y-\delta \Delta_Y$ is pseudo-effective, we can choose a positively curved metric $e^{-\varphi_\delta}$ on $F_\delta$.
Let $\varphi_{\pi^*\Delta_Y}$ (resp. $\varphi_R$) be the weight of the canonical singular metric on $\mathcal{O}_X (\pi^*\Delta_Y)$ (resp. $\mathcal{O}_X (R)$).
The minimality of the singularities of $\varphi_{F,\rm min}$ implies that
\begin{equation}\label{eq91}
 \varphi_\delta+\delta\varphi_{\pi^*\Delta_Y}\leq\varphi_{F,\rm min} +C
\end{equation}
holds for some constant $C$.
\smallskip

Define a metric $h_S$ on $\cO(S)$ by the weight
\begin{equation}\label{eq102}
 \varphi_S:=\frac{\varphi_{F,\rm min}-\varphi_\delta}{\delta \ep_0}
                -\frac{\varphi_R}{\ep_0}.
\end{equation}
As a consequence of \eqref{eq91}, we have
\begin{equation}\label{eq103}
 \psi:=\log|f_S|_{h_S} ^2
 =\frac{\varphi_\delta+\delta\varphi_{\pi^*\Delta_Y}-\varphi_{F,\rm min}}{\delta \ep_0} \leq C_1 \qquad\text{ on } X .
\end{equation}
Here $f_S \in H^0 (X, \mathcal{O}_X (S))$ denotes the canonical section of $\mathcal{O}_X (S)$, namely $\rm Div (f_S) =[S]$.
\medskip

\step{2}{A smooth extension}
Let $v\in H^0(S,mF|_S)$, and set
$$u=s_E^m \cdot v.$$
Choose a Stein cover $\{\Omega_i\}$ of $X$, and let $V_i$ be a holomorphic extension of $v$ to $\Omega_i$. Set $U_i=s_E^mV_i$.
Let $\sum \theta_i \equiv 1$ be a partition of unity. We can define a smooth global extension
$\wt U :=\sum_i \theta_i U_i$ of $u$. We have $\wt U=s_E^m \cdot \wt V$, where $\wt V =\sum_i \theta_i V_i$.

We follow the extension argument of \cite{CDM,Dem16}. Write
$$
 m(K_X+S+B)=K_X+S+L,\qquad
 L:=B+(m-1)(F+E),
$$
and equip $L$ with the metric of weight
\begin{equation}\label{eq92}
 \varphi_L=\varphi_B+(m-1)\varphi_{F,\rm min} +(m-1)\varphi_E.
\end{equation}
Notice that locally we have the equality
\begin{equation}\label{eq93}
 \frac{|U_i-U_j|^2}{|f_S|^2}e^{-\varphi_L}
 =\frac{|V_i-V_j|^2}{|f_S|^2}|s_E|^2
        e^{-\varphi_B-(m-1)\varphi_{F,\rm min}} .
\end{equation}
By construction, the quotient $\displaystyle \frac{V_i-V_j}{f_S}$ is holomorphic.
Moreover, by \eqref{nulelong},
$e^{-p\varphi_{F,\rm min}}$ is locally integrable for every $p>0$.
Since $B$ is klt, H\"older's inequality gives
$$\int_X  e^{-\varphi_B-(m-1)\varphi_{F,\rm min}} <+\infty.$$
Thus the quotient $\displaystyle \frac{|\db\wt U|^2}{|f_S|^2}$ is smooth, and
\begin{equation}\label{eq94}
 \int_X\frac{|\db\wt U|^2}{|f_S|^2}e^{-\varphi_L}<\infty
\end{equation}
follows.

\medskip

\step{3}{Simultaneous regularization}
By construction, we have
$$i\Theta_{\varphi_L} (L) \geq 0 \qquad\text{ and } \qquad i\Theta_{\varphi_L} (L) \geq \delta\ep_0 (m-1) \cdot i\Theta_{\varphi_S} (\mathcal{O} (S)).$$
Together with \eqref{eq103}, these inequalities give the curvature and norm bounds required for the Ohsawa--Takegoshi-type extension theorems of \cite[Theorem~2.1]{CDM} and \cite[Section~4]{DHP}, apart from the regularity of the metrics.
We use Demailly's approximation of positive currents by Bergman kernels; see \cite[Chapter 13]{Dem12}.

Fix one sufficiently positive ample line bundle $A_X$ on $X$, with a smooth, positively curved metric. For all sufficiently divisible $k$, we consider the  spaces
\[
\mathcal H_{F,k}=H^0\bigl(X,(kF+A_X)\otimes\mathcal{I}(k\varphi_{F,\min})\bigr),
\quad
\mathcal H_{\delta,k}=H^0\bigl(X,(kF_\delta+A_X)\otimes \mathcal{I}(k\varphi_\delta)\bigr),
\]
with the $L^2$ norms induced by the indicated weights and a fixed smooth volume form. Let $B_{F,k},B_{\delta,k}$ be the corresponding Bergman kernels, and define the local weights on the original rational bundles by
\[
\varphi_{F,k}=\frac1k(\log B_{F,k}-\varphi_{A_X}),\qquad
\varphi_{\delta,k}=\frac1k(\log B_{\delta,k}-\varphi_{A_X}).
\]
Choose $A_X$ sufficiently ample for the standard global Bergman approximation estimates and the global generation of the sheaves defining these systems. Such a choice can be made uniformly in $k$, for example by $L^2$ estimates.

By \eqref{nulelong}, $\mathcal{I}(k\varphi_{F,\min})=\cO_X$. It follows immediately that $\varphi_{F,k}$ is smooth. The regularizations converge in $L^1$ and, after selecting a common subsequence, they converge almost everywhere. By construction,  we have
\[
\sqrt{-1}\Theta_{\varphi_{F,k}}(F)\geq-\frac{C_0}{k}\omega,
\qquad
\sqrt{-1}\Theta_{\varphi_{\delta,k}}(F_\delta)\geq-\frac{C_0}{k}\omega.
\]
Their normalizations can be chosen so that
\[
\varphi_{F,k}\geq\varphi_{F,\min}-C_0,
\qquad \sup_X\varphi_{F,k}\leq C_0,
\qquad \sup_X\varphi_{\delta,k}\leq C_0.
\]
Here the upper bounds are understood using fixed smooth reference metrics, or on a fixed finite relatively compact trivializing cover.
\smallskip

Multiplication by $s_{\pi^*\Delta_Y}^{k\delta}$ maps $\mathcal H_{\delta,k}$ into $\mathcal H_{F,k}$. Moreover, the norm of this operator is smaller than $e^{kC}$ for some positive constant $C> 0$ independent of $k$, by \eqref{eq91}. The extremal characterization of the Bergman kernels therefore yields
\begin{equation}\label{eq97}
\delta\varphi_{\pi^*\Delta_Y}+\varphi_{\delta,k}\leq\varphi_{F,k}+C,
\end{equation}
uniformly with respect to $k$.

Define the regularized weights and functions by
\begin{align}
 \varphi_{S,k}
 &=\frac{\varphi_{F,k}-\varphi_{\delta,k}}{\delta \ep_0}
       -\frac{\varphi_R}{\ep_0},\label{eq:regularized-S}\\
 \psi_k
 &=\log|f_S|^2-\varphi_{S,k},\label{eq:regularized-psi}\\
 \varphi_{L,k}
 &=\varphi_B+(m-1)\varphi_E+(m-1)\varphi_{F,k}.
 \label{eq:regularized-L}
\end{align}
As in \eqref{eq103}, inequality \eqref{eq97} gives
\begin{equation}\label{eq977}
\psi_k\leq C
\end{equation}
for some constant $C$ independent of $k$.
The curvature conditions become
\begin{equation}\label{eq:curv1}
  \sqrt{-1}\Theta_{\varphi_{L,k}} (L)\geq-\frac{C_1}{k}\omega,
  \end{equation}
  \begin{equation}\label{eq:curv2}
  \sqrt{-1}\Theta_{\varphi_{L,k}} (L)- \delta\ep_0 (m-1) \cdot  \sqrt{-1}\Theta_{\varphi_{S,k}} (S)
 \geq-\frac{C_1}{k}\omega.
 \end{equation}
for some constant $C_1$ independent of $k$.

Moreover, since $\varphi_{F,k}$ is smooth,
$\psi_k = -\infty$ on $S$. By the construction of $\psi_k$, we have
\begin{equation}\label{eq98}
 e^{-\varphi_{L,k}-\varphi_{S,k}-\psi_k} =\frac{e^{-\varphi_{L,k}}}{|f_S|^2}.
\end{equation}
Thus the auxiliary weight $\varphi_{S,k}$ cancels from the norm used in the approximate equation.
\smallskip

\medskip

\step{4}{Uniform estimates} For convenience, set
$$ I_{t,k}: =\int_{\{t<\psi_k<t+1\}} \frac{|\wt U|^2}{|f_S|^2}e^{-\varphi_{L,k}}.
$$
On the integration domain, as $e^{-\psi_k}\leq e^{-t}$, we have
\begin{equation}\label{eq99}
 I_{t,k}\leq e^{-t}\int_X
           |\wt U|^2e^{-\varphi_{L,k}-\varphi_{S,k}}.
\end{equation}
Since $\varphi_{F,\rm min}$ has zero Lelong numbers everywhere and $B$ is klt, H\"older's inequality gives
$$
 \int_X |\wt U|^2e^{-\varphi_{L,k}-\varphi_{S,k}} = \int_X |\wt V|^2|s_E|^2|s_R|^{2/\ep_0}e^{-\varphi_B}
 \exp\left[-\left(m-1+\frac1{\delta \ep_0}\right)\varphi_{F,k}
                    +\frac{\varphi_{\delta,k}}{\delta \ep_0}\right]
$$
$$\leq C_1 \int_X e^{-\varphi_B-\left(m-1+\frac1{\delta \ep_0}\right)\varphi_{F,\min}} < +\infty .$$

\noindent Therefore we have
\begin{equation}\label{eq101}
I_{t,k}\leq C \cdot e^{-t}
\end{equation}
for some constant $C$ independent of $t$ and $k$.

\medskip

\step{5}{Approximate solutions}
This is the crucial step.
We adapt the argument of \cite[Theorem~2.1, Step~2]{CDM} to the total weight
$\varphi_{L,k}+\varphi_{S,k}+\psi_k$, recording the approximate-solution estimate in the form needed here.
\medskip

\noindent We denote by $a:=\delta\varepsilon_0(m-1)>0$ the constant which appears in \eqref{eq:curv2}. After a common constant
normalization of the metric on $S$, we can assume that $\psi_k\leq0$.
The regularized curvature inequalities needed here are \eqref{eq:curv1} and \eqref{eq:curv2}.
Note that the weight
$\varphi_{S,k}$ is not required to be quasi-psh.

Next we make a specific choice of a cutoff function, depending on a fixed
parameter $0<\gamma\leq a$. For $t\leq-2$ and $x\leq0$, define
\[
\chi_t(x):=t+\log(1+e^{x-t})-\log2,
\]
and put
\[
\eta_{t,k}=1-\gamma\chi_t(\psi_k),\qquad
\lambda_{t,k}=2\gamma
\frac{\chi_t'(\psi_k)^2}{\chi_t''(\psi_k)}
=2\gamma e^{\psi_k-t}.
\]
These are positive smooth functions in the complement of the analytic singularities.
Since
\[
t-\log2\leq\chi_t(x)\leq0,\quad 0<\chi_t'(x)<1,
\quad
\chi_t''(x)=\frac{e^{x-t}}{(1+e^{x-t})^2},
\]
we have
\[
1\leq\eta_{t,k}\leq 1+\gamma(|t|+\log2),\qquad
0<\lambda_{t,k}\leq2\gamma e^{-t}.
\]
In particular, both functions are bounded for fixed $t$, uniformly with respect to $k$.
When restricted to the set $\{t<\psi_k<t+1\}$, the inequality
\[
\chi_t''(\psi_k)\geq b_0:=\frac{e}{(1+e)^2}>0
\]
holds.

Let $Z_k$ be an ample divisor of sufficiently high degree containing $S$, the supports of $B,E$, and the analytic
singularities of the regularized weights. Its complement admits a
complete K\"ahler metric. We use the original restricted $\omega$
in the norms, which is allowed in the twisted estimate. On this
complement the total metric has weight
\[
H_k=\varphi_{L,k}+\varphi_{S,k}+\psi_k
=\varphi_{L,k}+\log|f_S|^2,
\]
so its curvature equals $\sqrt{-1}\Theta_{\varphi_{L,k}} (L)$. A direct calculation gives
\begin{align*}
R_{t,k}
&:=\eta_{t,k}\sqrt{-1}\Theta_{H_k}-\sqrt{-1}\partial\bar\partial\eta_{t,k}
-\lambda_{t,k}^{-1}\sqrt{-1}\partial\eta_{t,k}\wedge\bar\partial\eta_{t,k}\\
&=\eta_{t,k}\left(\sqrt{-1}\Theta_{\varphi_{L,k}} (L)
-\frac{\gamma\chi_t'(\psi_k)}{\eta_{t,k}}\sqrt{-1}\Theta_{\varphi_{S,k}} (S)\right)
+\frac{\gamma\chi_t''(\psi_k)}2
\sqrt{-1}\partial\psi_k\wedge\bar\partial\psi_k.
\end{align*}
The ratio in the parentheses belongs to $[0,a]$. By interpolating the
two curvature conditions \eqref{eq:curv1} and \eqref{eq:curv2}, we obtain
\[
R_{t,k}\geq-\frac{C(1+|t|)}{k}\omega
+\frac{\gamma\chi_t''(\psi_k)}2
i\partial\psi_k\wedge\bar\partial\psi_k,
\]
where $C$ is independent of $k,t$.

Let $\theta\colon\mathbb R\to[0,1]$ be smooth, equal to $1$ on $(-\infty,0]$ and to $0$ on $[1,\infty)$. We write
\[
v_{t,k}=\bar\partial(\theta(\psi_k-t)\widetilde U)
=v^{(1)}_{t,k}+v^{(2)}_{t,k},
\]
\[
v^{(1)}_{t,k}=\theta'(\psi_k-t)\bar\partial\psi_k\wedge\widetilde U,
\qquad
v^{(2)}_{t,k}=\theta(\psi_k-t)\bar\partial\widetilde U,
\]
and use the approximate twisted estimate with error
$\displaystyle \epsilon_{t,k}=\frac{2C(1+ |t|)}{k}$. Note that the positive term in the lower bound of $R_{t, k}$ controls
$v^{(1)}_{t,k}$ uniformly on the support of $\theta'(\psi_k-t)$.
By applying the $L^2$ estimate \cite[Prop 1.3]{CDM}, we can find sections $u_{t,k}$ and $w_{t,k}$ satisfying
\begin{equation}\label{eq104}
\bar\partial u_{t,k}+w_{t,k}=v_{t,k}
\end{equation}
and
\begin{equation}\label{eq:expliciterror}
\int_X\frac{|w_{t,k}|^2}{|f_S|^2}e^{-\varphi_{L,k}}
\leq\frac{C(1+ |t|)}{k}I_{t,k}
+C\int_{\{\psi_k<t+1\}}
\frac{|\bar\partial\widetilde U|^2}{|f_S|^2}e^{-\varphi_L}.
\end{equation}
Note that \eqref{eq:expliciterror} is a consequence of \cite[(18)]{CDM}, \eqref{eq98}, combined with the fact that
$e^{-\varphi_{L,k}}\leq C e^{-\varphi_L}$.

The estimate also gives
\begin{equation}\label{eq105}
\int_X(\eta_{t,k}+\lambda_{t,k})^{-1}
\frac{|u_{t,k}|^2}{|f_S|^2}e^{-\varphi_{L,k}}<\infty.
\end{equation}
In particular $\displaystyle \frac{u_{t,k}}{f_S}$ and $\displaystyle \frac{w_{t,k}}{f_S}$ are locally $L^2$ for smooth metrics.
Initially, \eqref{eq104} holds on $X\setminus Z_k$. The $L^2$ condition extends it across $Z_k$, so it holds on all of $X$ in the sense of distributions.

\medskip

\noindent We now choose $k$ as a function of $t$ so that the error $w_{t,k}$ tends to zero as $t\to-\infty$. Define the finite measure
\[
d\mu:=\frac{|\bar\partial\widetilde U|^2}{|f_S|^2}
e^{-\varphi_L}dV_\omega
\]
cf. \eqref{eq94}. Along an almost-everywhere convergent
regularization sequence, $\mu(\{|\psi_k-\psi|>1\})\to0$.
Choose $t_j=-j$ for $j\geq2$. For each $j$, choose a
large and sufficiently divisible $k_j$ so that
\[
\mu(\{|\psi_{k_j}-\psi|>1\})\leq\frac1j,
\qquad
\frac{C(1+ |t_j|)e^j}{k_j}\leq\frac1j.
\]
For simplicity, write
$$w_j :=w_{t_j, k_j} \qquad\text{and} \qquad u_j:= u_{t_j, k_j} .$$
Since $I_{t,k}\leq Ce^{-t}$ by the previous step, and
\[
\{\psi_{k_j}<t_j+1\}\subset
\{\psi<t_j+2\}\cup\{|\psi_{k_j}-\psi|>1\},
\]
equation \eqref{eq:expliciterror} implies
\[
\int_X\frac{|w_j|^2}{|f_S|^2}e^{-\varphi_{L,k_j}}
\leq\frac Cj+C\mu(\{\psi<t_j+2\})\longrightarrow0.
\]
Here $\psi$ is finite almost everywhere, so the last tail tends to
zero. The uniform upper bound for $\varphi_{L,k_j}$ then gives
\begin{equation}\label{eq1071}
\lim_{j\to +\infty}\int_X \Big|\frac{w_j}{f_S}\Big|^2 dV_{\omega} =0
\end{equation}
where the norm $|\cdot |$ is with respect to a smooth metric.
\smallskip

\begin{rmk}
Our choice of truncation function $\eta_{t,k}$ gives an explicit negative lower bound for the twisted curvature term $R_{t,k}$. The argument would also work with an arbitrary finite function $C(t)$ in place of $C(1+|t|)$ in \eqref{eq:expliciterror}. We retain the explicit bound, as it may be useful in other contexts.
\end{rmk}

\medskip

\step{6}{Vanishing of the obstruction class}
Consider the smooth $L$-valued $(n,1)$-form
$$\alpha :=\frac{\db\wt U}{f_S} $$
which represents the extension obstruction.  From the relation
\eqref{eq104} and the fact that $\psi_k =-\infty$ on $S$, we obtain
\begin{equation}\label{eq108}
 \alpha-\frac{w_j}{f_S}
 =\db\left(
 \frac{(1-\theta(\psi_{k_j}-t_j))\wt U+u_j}{f_S}
 \right) \qquad\text{ on } X
\end{equation}
in the sense of currents.
Initially, \eqref{eq108} is defined outside $S$. Since $w_j/f_S$ and $u_j/f_S$ are locally $L^2$, the standard removability lemma for $L^2$ forms extends the identity across $S$.
Therefore we have
$$ [\alpha]=\Big[\frac{w_j}{f_S}\Big] \in H^1 (X, K_X+L).$$
With respect to a fixed smooth metric, the harmonic projections of the $\dbar$-closed $L^2$ forms $w_j/f_S$ tend to zero by \eqref{eq1071}. Since each is the harmonic representative of $[\alpha]$, we obtain
$$
[\alpha]=0\in H^1(X, K_X+L).
 $$
 Therefore there exists a smooth $L$-valued $(n,0)$-form $\beta$ such that
$\db\beta=\alpha$. Then
$$
 \sigma :=\wt U-f_S \cdot \beta \in H^0 (X, K_X +S+L) = H^0 (X, m F+ mE) .
$$

\noindent Finally, set $\displaystyle \widehat\sigma:=\frac{\sigma}{s_E^m}$. Since $E$ is exceptional, the section $\widehat\sigma$ is holomorphic. Moreover, since $S$ and $E$ have no common component, $\widehat\sigma|_S=v$. This completes the proof.
\end{proof}

\begin{rmk}
The extension argument above only uses pseudo-effectivity of
$F_Y-\delta\Delta_Y$ and the zero-Lelong-number residual metric
provided by Theorem~\ref{nu0}. It does not require the existence of the signed divisor
representing $F_Y$ with support in the boundary.
\end{rmk}
\medskip

\noindent We conclude with the following consequence of Theorem~\ref{extensionthm} and \cite[Lemma~3.2]{OpenAI2}. For the reader's convenience, we include an almost self-contained proof.

\begin{cor}[\cite{OpenAI2}, Theorem~3.1]\label{coropenAI}
Let $(Y,\Delta_Y)$ satisfy the hypotheses of Theorem~\ref{extensionthm}, with $\Delta_Y$ reduced. Assume that
\[
K_Y+ \Delta_Y\sim_{\mathbb Q}G:=\sum_i a_iD_i,\qquad a_i\in\mathbb Q,
\qquad \Delta_Y=\sum_iD_i,
\]
and that $(K_Y+\Delta_Y)|_{\Delta_Y}$ is semiample as a $\bbQ$-line bundle on the whole reduced scheme $\Delta_Y$. Then the equality
\[
\kappa(Y,K_Y+ \Delta_Y)=\nu(Y,K_Y+ \Delta_Y)
\]
holds.
\end{cor}

\begin{proof}
We begin by introducing some notation. On the chosen log resolution
$\pi:X\to Y$, set
\[
F=\pi^\star(K_Y+\Delta_Y),\qquad
\Xi:=(\supp \pi^\star\Delta_Y)_{\mathrm{red}},\qquad
G_X=\pi^\star G.
\]
The divisor $\Xi$ has simple normal crossings, contains $S$, and contains the support of $G_X$. We can write
\[
G_X=\sum_\alpha\beta_\alpha\Xi_\alpha,
\]
where the $\Xi_\alpha$ range over all irreducible components of $\Xi$.
Let $\Xi_+$, $\Xi_-$, and $\Xi_0$ denote the reduced unions of the components of $\Xi$ corresponding to positive, negative, and zero coefficients $\beta_\alpha$, respectively.

Since $\pi^\star\Delta_Y$ has strictly positive coefficients along every component of $\Xi$, we can choose $c>0$ such that
$\delta\pi^\star\Delta_Y-c\Xi\geq0$, where $\delta$ is as in Theorem~\ref{extensionthm}. Thus,
\[
F-c\Xi=(F-\delta\pi^\star\Delta_Y)
       +(\delta\pi^\star\Delta_Y-c\Xi)
\]
is pseudo-effective.
\smallskip

\noindent
Choose a sufficiently divisible integer $q>0$ such that
$q(K_Y+\Delta_Y)|_{\Delta_Y}$ is generated by its global sections.
Let $\phi:\Delta_Y\to\mathbb P^N$ be the corresponding morphism.
Since $\pi(\Xi)=\Delta_Y$, we can define
$\Phi:=\phi\circ\pi|_\Xi$. We then have
\begin{equation}\label{oai1}
\Phi:\Xi\longrightarrow\mathbb P^N,
\qquad
qF|_\Xi=\Phi^\star\mathcal O_{\mathbb P^N}(1).
\end{equation}
In particular, $F|_\Xi$ is semiample.
\smallskip

\noindent
We divide the remainder of the argument into several steps, with the aim of making the ideas in \cite{OpenAI2} and their interaction with Theorem~\ref{extensionthm} more transparent.

\step{1}{Properties of $\Xi$}
We recall the numerical argument from \cite[Lemma~3.2]{OpenAI2}.
Fix an ample divisor $H$ on $X$ and set $v:=\nu(F)$.
By the nefness of $F$ and the pseudo-effectivity of $F-c\Xi$, we obtain
\[
0\leq (F-c\Xi)F^vH^{n-v-1}
=-c\sum_\alpha \Xi_\alpha F^vH^{n-v-1}.
\]
It follows that
$\Xi_\alpha\cdot F^v\cdot H^{n-v-1}=0$
for every $\alpha$.
Together with \eqref{oai1}, this implies that
\[
\dim\Phi(\Xi_\alpha)\leq v-1
\quad\text{for every }\alpha.
\]
On the other hand,
\[
0<F^vH^{n-v}
=\sum_\alpha\beta_\alpha\Xi_\alpha F^{v-1}H^{n-v}.
\]
Hence, there exists a component $\Xi_\alpha\subset\Xi_+$
such that $\dim\Phi(\Xi_\alpha)=v-1$.
\smallskip

We next show that
\begin{equation}\label{eq:smallintersection}
\dim\Phi(J)<v-1,
\end{equation}
where $J:=\Xi_+\cap(\Xi_-\cup\Xi_0)$.
Consider the quadratic form
\[
\mathcal Q(U,V)=U \cdot V\cdot F^{v-1} \cdot H^{n-v-1}.
\]
By the preceding discussion, we have
\[
\mathcal Q(F,F)=0,\qquad
\mathcal Q(F,\Xi_\alpha)=0,\qquad
\mathcal Q(F,H)=F^vH^{n-v}>0.
\]
It follows from the mixed Hodge index theorem that the matrix
$M:=(\mathcal Q(\Xi_\alpha,\Xi_\gamma))_{\alpha,\gamma}$
is negative semidefinite.\footnote{Indeed, let
$X:=\sum_\alpha a_\alpha \Xi_\alpha$ be an arbitrary linear combination.
Since $\mathcal Q(F,H)>0$, there exists a constant $c$ such that
$\mathcal Q(X-cF,H)=0$.
By the mixed Hodge index theorem, we have
$\mathcal Q(X-cF,X-cF)\leq 0$.
Together with $\mathcal Q(X,F)=0$ and $\mathcal Q(F,F)=0$,
this implies $\mathcal Q(X,X)\leq 0$. Therefore $M$ is negative semidefinite.}

Write $G_X=(G_X)_++(G_X)_-$, where $(G_X)_+$
(resp. $(G_X)_-$) is the sum of the components of $G_X$
with positive (resp. negative) coefficients.
Let $\beta=(\beta_\alpha)_\alpha$ be the vector of
coefficients of $G_X$.
Since $\mathcal Q(F,F)=0$, we have
\begin{equation}\label{adddbeta}
M\beta=0.
\end{equation}
Similarly, let $\beta^+$ and $\beta^-$ denote the vectors
of coefficients of $(G_X)_+$ and $(G_X)_-$, respectively.
Then $\beta=\beta^++\beta^-$.

Since $\mathcal Q(\Xi_\alpha,\Xi_\gamma)\geq 0$
for $\alpha\neq\gamma$, and the nonzero coefficients of $\beta^-$ are negative,
we have
\[
(\beta^+)^tM\beta^-\leq 0.
\]
Recall that by \eqref{adddbeta}, we have
$$0= (\beta^+)^tM \beta = (\beta^+)^tM \beta^+ + (\beta^+)^tM  \beta^- .$$
Since both terms on the right are nonpositive, we obtain
\[
M\beta^+=0.
\]
Consequently,
\[
\mathcal Q((G_X)_+,\Xi_\alpha)=0
\quad\text{for every }\Xi_\alpha\subset\Xi_-\cup\Xi_0.
\]
This yields \eqref{eq:smallintersection}.

\step{2}{Sections obtained by extension}
We first show that there exists a component $\Xi_\star$ of $\Xi_+$
such that $\dim\Phi(\Xi_\star)=v-1$ and $\Xi_\star$ is the
strict transform of some $D_i$. In particular, $\Xi_\star$ is a
component of $S$.

Indeed, by Step~1, there exists a component $\Xi_\gamma$
of $\Xi_+$ such that $\dim\Phi(\Xi_\gamma)=v-1$. Since
\[
\beta_\gamma
=\sum_i a_i\,\ord_{\Xi_\gamma}(\pi^\star D_i)>0,
\]
there exists an index $i$ such that $a_i>0$ and
$\ord_{\Xi_\gamma}(\pi^\star D_i)>0$. Consequently,
\[
\Phi(\Xi_\gamma)\subset\phi(D_i)
=\Phi(\widehat D_i),
\]
where $\widehat D_i$ denotes the strict transform of $D_i$.
Set $\Xi_\star:=\widehat D_i$. Then $\Xi_\star$ has the desired
properties: indeed, $\dim\Phi(\Xi_\star)=v-1$, and its
coefficient in $G_X$ is $a_i>0$.

\medskip

By Step~1, we have $\dim\Phi(J)<v-1$, whereas
$\dim\Phi(\Xi_\star)=v-1$.
Together with the relation
$qF|_\Xi=\Phi^\star\mathcal O_{\mathbb P^N}(1)$,
this implies that, for a sufficiently divisible integer $m_0>0$,
we can find sections
\[
v_0,\ldots,v_N\in H^0(\Xi,m_0F|_\Xi)
\]
such that $v_0|_{\Xi_\star}\not\equiv 0$,
$v_i|_{\Xi_-\cup\Xi_0}=0$ for every $i$, and the rational
map defined by
$[v_0,\ldots,v_N]|_{\Xi_\star}$
has image of dimension $v-1$. Here we use the simple normal crossings property of $\Xi$: the scheme-theoretic intersection $\Xi_+\cap(\Xi_-\cup\Xi_0)$ is reduced. Sections vanishing on this intersection therefore extend by zero to $\Xi_-\cup\Xi_0$.
\smallskip

\noindent By Theorem~\ref{extensionthm}, each $v_i|_S$ extends to $X$, giving
\[
\sigma_i\in H^0(X,m_0F),\qquad \sigma_i|_S=v_i|_S.
\]

\medskip

We claim that $\sigma_0$ \emph{vanishes along every component of $\Xi_-$, including those which are not part of the support of $S$}.

To begin with, let $D_i \subset \Delta_Y$ be a component with $a_i<0$, and let $\widehat D_i$ denote its strict transform. Then $\widehat D_i \subset \Xi_-$. Since $\Delta_Y$ is reduced, we also have $\widehat D_i \subset S$. By the construction of $\sigma_0$, it follows that
\begin{equation}\label{addpoint}
\sigma_0|_{\widehat D_i}=0
\end{equation}
for every $i$ such that $a_i<0$.

Now let $\Xi_\alpha$ be an irreducible component of $\Xi_-$, and let $\beta_\alpha<0$ denote its coefficient. We now prove that $\sigma_0|_{\Xi_\alpha}=0$. We distinguish two cases.

\smallskip
\noindent\textbf{Case 1:} $\Xi_\alpha$ is the strict transform $\widehat D_\alpha$ of some $D_\alpha$.  In this case, $a_\alpha=\beta_\alpha<0$. By \eqref{addpoint}, we obtain
$\sigma_0|_{\Xi_\alpha}=0$.

\smallskip
\noindent\textbf{Case 2:} $\Xi_\alpha$ is $\pi$-exceptional.
Recall that
\begin{equation}\label{eq:oai2}
\sigma_0=\pi^\star\sigma_{0,Y}.
\end{equation}
Together with \eqref{addpoint}, this implies that
\begin{equation}\label{eq:oai22}
\sigma_{0,Y}|_{D_i}=0
\end{equation}
for every $i$ such that $a_i<0$.

On the other hand, since
\[
\sum_i a_i\,\ord_{\Xi_\alpha}(\pi^\star D_i)
=\beta_\alpha<0,
\]
there exists at least one index $i$ such that $a_i<0$ and
$\ord_{\Xi_\alpha}(\pi^\star D_i)>0$.
Consequently, $\pi(\Xi_\alpha)\subset D_i$ for some $i$ with $a_i<0$.
Combining this with \eqref{eq:oai2} and \eqref{eq:oai22}, we conclude that
$\sigma_0|_{\Xi_\alpha}=0$. The claim is thus proved.

\smallskip

\step{3}{An additional section vanishing on $\Xi_\star$}
Since $F\sim_{\mathbb Q} G_X$, for a sufficiently divisible positive integer $q_0$, there exists a rational section $\tau$ of $q_0F$ whose associated divisor is $q_0G_X$. Thus, the poles of $\tau$ are supported on $\Xi_-$, while $\tau$ vanishes along $\Xi_\star$.
As shown in Step~2, the section $\sigma_0$ vanishes along $\Xi_-$. Hence, for a sufficiently large integer $N_0$, the section
$$
\rho:=\tau^{m_0}\sigma_0^{N_0}\in H^0(X,m_0(q_0+N_0)F)
$$
is a nonzero holomorphic section that vanishes along $\Xi_\star$.
\medskip

\step{4}{Counting} Introduce
\[
 d_0:=m_0(q_0+N_0),\qquad
 u_i:=\sigma_0^{q_0+N_0-1}\sigma_i\in H^0(X,d_0F).
\]
The restriction of $u_0$ to $\Xi_\star$ is not identically zero. The restricted ratios
$(u_i/u_0)|_{\Xi_\star}=(\sigma_i/\sigma_0)|_{\Xi_\star}$ generate a field of transcendence degree $v-1$.
After relabeling, choose $u_1,\ldots,u_{v-1}$ so that
\[
 \left.\frac{u_1}{u_0}\right|_{\Xi_\star},\ldots,
 \left.\frac{u_{v-1}}{u_0}\right|_{\Xi_\star}
\]
are algebraically independent over $\mathbb C$.

Recall that $\rho\in H^0(X,d_0F)$ is nonzero and vanishes along $\Xi_\star$.
For every $k\ge1$, every
$0\le j\le k$, and every $\alpha=(\alpha_1,\ldots,\alpha_{v-1})\in\bbN^{v-1}$
with $|\alpha|\le k-j$, define
\[
 s_{j,\alpha}:
 =\rho^j u_0^{\,k-j-|\alpha|}
       u_1^{\alpha_1}\cdots u_{v-1}^{\alpha_{v-1}}
 \in H^0(X,kd_0F).
\]
The vanishing of $\rho$ along $\Xi_\star$, together with the algebraic independence of the restricted ratios, shows that these sections are linearly independent. Counting them gives $\kappa(F)=v$, proving the corollary.

\end{proof}

\section{Appendix}

\noindent We give an $L^2$ proof of Lemma~\ref{small}, beginning with a more general statement in which the base is only assumed to be normal.

Let $Y$ be a normal projective variety, and let $\pi:X\to Y$ be a projective resolution. Let $F_Y$ be a nef $\mathbb Q$-Cartier divisor, and set $F:=\pi^\star F_Y$. Assume that
\begin{equation}\label{eq:adjoint}
K_X+S+B=F+E,
\end{equation}
where $S$ is reduced, $B,E\geq0$ are rational divisors, the coefficients of $B$ belong to $[0,1)$, and $E$ is $\pi$-exceptional. Assume that $S,B,E$ have pairwise disjoint sets of irreducible components and that their union with the exceptional locus has simple normal crossings. We fix $x_0\in X$, and put $y_0=\pi(x_0)$. All sufficiently divisible integers below are understood to clear the relevant denominators.
\smallskip

\noindent We start with the following statement.

\begin{lem}\label{Appnef}
There exist an ample line bundle $A_Y$ on $Y$ and an integer $C_0\geq0$ such that, for every sufficiently divisible integer $m\geq2$, there is a nonzero section
\begin{equation}\label{upperboundord}
\sigma_m\in H^0(X,mF+\pi^*A_Y)
\quad\text{with}\quad
\ord_{x_0}\sigma_m\leq C_0.
\end{equation}
\end{lem}

\begin{proof}
 We may assume that $x_0\in\supp S$; otherwise, the same argument applies with minor simplifications. Set $y_0:=\pi(x_0)$, and choose a small Stein neighborhood $V$ of $y_0$ in $Y$ on which the bundles on $Y$ under consideration are trivial.

The proper image $\pi(S)$ is a proper closed subset of $Y$. We can therefore find a nonzero holomorphic function $g$ on $V$, vanishing on $\pi(S)\cap V$. Since $S$ is reduced, $\pi^*g$ vanishes along $S\cap\pi^{-1}(V)$.

Fix a very ample line bundle $H$ on $Y$. We choose a smooth cutoff $\rho$ on $Y$, with compact support in $V$ and equal to one near $y_0$. Here and below smooth functions and metrics on the singular space $Y$ are understood by restriction from local ambient embeddings. If $\omega_H$ is the restricted Fubini--Study form of the embedding defined by $H$, then the inequality
\begin{equation}\label{eq:cutoffbound}
\sqrt{-1}\partial \rho\wedge\bar\partial\rho
\leq C_\rho\omega_H
\end{equation}
holds for a fixed constant $C_\rho$.

Choose a function $\psi_{y_0}$ with a logarithmic pole at $y_0$ satisfying
\begin{equation}\label{eq:polecurvature}
\ddbar\psi_{y_0}\geq-\omega_H.
\end{equation}
Such a function can be constructed explicitly on the ambient projective space and then restricted to $Y$.
\smallskip

Choose integers
\[b>C_\rho+1, \qquad a>\lambda+ b\]
where
\[\lambda:=n+\ord_{x_0}(\pi^\star g)+\ord_{x_0}\lceil E\rceil,\]
and set $A_Y=aH$. For each $m$, let $e_m$ be a frame of $mF_Y+A_Y$ on $V$, and define
\begin{equation}\label{eq:localsection}
s_m=s_{mE}\,(\pi^\star g)\,\pi^\star e_m
\in H^0\bigl(\pi^{-1}V,mF+mE+\pi^*A_Y\bigr).
\end{equation}
By construction,
\[
\Div(s_m)\geq mE+S,
\qquad
\ord_{x_0}\frac{s_m}{s_{mE}}=\ord_{x_0}(\pi^\star g).
\]
Through \eqref{eq:adjoint}, this is a local section of $K_X+L_m$, where
\[
L_m=S+B+(m-1)E+(m-1)F+\pi^*A_Y.
\]

Next, the $\bbQ$-divisor
\[
(m-1)F_Y+(a-\lambda-b)H
\]
is ample. Choose a smooth positively curved metric denoted by $e^{-\tau_m}$ on it, whose pullback is smooth on $X$. We then endow $L_m$ with the singular metric of weight
\begin{equation}\label{eq:metric}
\Phi_m=\varphi_S+\varphi_B+(m-1)\varphi_E
 +\pi^\star\tau_{m}+(\lambda+b)\pi^*\varphi_H
 +\lambda\psi_{y_0}\circ\pi.
\end{equation}
By \eqref{eq:polecurvature},
\begin{equation}\label{eq:positive}
\sqrt{-1}\Theta_{\Phi_m}(L_m)
\geq[S]+[B]+(m-1)[E]+b\pi^*\omega_H.
\end{equation}
The smooth pullback contribution is semipositive on $X$, and is positive definite off the exceptional locus.

Consider the $L_m$-valued $(n,1)$-form
\[
v_m=\bar\partial\bigl((\rho\circ\pi)s_m\bigr)
=\bar\partial(\rho\circ\pi)\wedge s_m.
\]
On the support of $\bar\partial(\rho\circ\pi)$, the log-pole function $\psi_{y_0}\circ\pi$ is bounded, and the quotient $\pi^*g/f_S$ is holomorphic. We then have
\[
|s_m|^2e^{-\varphi_S-\varphi_B-(m-1)\varphi_E}
\asymp
\left|\frac{\pi^*g}{f_S}\right|^2
|s_E|^2e^{-\varphi_B},
\]
where smooth factors are suppressed. This is locally integrable because $B$ has snc support and coefficients less than one. Together with \eqref{eq:cutoffbound} and \eqref{eq:positive}, this gives finiteness of the inverse-curvature integral for $v_m$.

The singular $L^2$ estimate therefore yields
\begin{equation}\label{eq:solution}
\bar\partial u_m=v_m,
\qquad
\int_X|u_m|^2e^{-\Phi_m}<\infty
\end{equation}
Thus
\[
\widetilde\sigma_m=u_m-(\rho\circ\pi)s_m
\in H^0(X,mF+mE+\pi^*A_Y)
\]
is holomorphic.

\noindent Since $Y$ is normal and $mE$ is effective and exceptional,
every such section is divisible by $s_{mE}$. Hence
\begin{equation}\label{eq:globalsection}
\sigma_m=\frac{\widetilde\sigma_m}{s_{mE}}
\in H^0(X,mF+\pi^*A_Y)
\end{equation}
is holomorphic.
\smallskip

The weighted integrability implies that
\[
\ord_{x_0}\frac{u_m}{s_{mE}}\geq \ord_{x_0}(\pi^\star g)+1.
\]
Indeed, $u_m$ is holomorphic near $x_0$, and the conclusion follows from the choice of $\Phi_m$ and \eqref{eq:solution}.
By contrast, $s_m/s_{mE}$ has order exactly $\displaystyle \ord_{x_0}(\pi^\star g)$. Since $\rho=1$ near $y_0$, these lowest-order terms cannot cancel in \eqref{eq:globalsection}. Therefore, setting $C_0:=\ord_{x_0}(\pi^\star g)$, we obtain
\[
\ord_{x_0}\sigma_m=C_0.
\]
In particular, $\sigma_m$ is nonzero, and the lemma follows.
\end{proof}
\medskip

\noindent Let $F$ be a pseudo-effective $\bbQ$-line bundle on a smooth projective manifold $X$.
Choose a positively curved metric $e^{-\varphi_F}$ on $F$ and a point $x_0$ such that $\nu(\varphi_F,x_0)=0$. For every $c>0$, there is a neighborhood $U_c$ of $x_0$ such that
\[x_0\in U_c, \qquad \exp(-c\varphi_F)\in L^1(U_c).\]
The same cutoff-and-correction argument, applied to the adjoint bundle $mF+A-K_X$ with the weight $m\varphi_F$ and a sufficiently positive fixed ample twist, yields a nonzero section of $mF+A$ for every sufficiently divisible $m$.

Finally, we proceed by interpolation as in the proof of Lemma~\ref{small}.

\bibliographystyle{alpha}

\end{document}